\documentclass[11pt,a4paper]{amsart} 
\usepackage{amssymb,bbm,bm}
\usepackage{epsfig}
\usepackage[all]{xy}

\newtheorem{theorem}{Theorem}[section] 
\newtheorem{lemma}[theorem]{Lemma}
\newtheorem{proposition}[theorem]{Proposition}

\newtheorem{corollary}[theorem]{Corollary}
\newtheorem{definition}[theorem]{Definition}
\newtheorem{remark}[theorem]{Remark}

\def\1{{\mathbbm{1}}}

\def\Aut{{\rm Aut}}
\def\be{\begin{equation}}
\def\build#1_#2^#3{\mathrel{\mathop{\kern 0pt#1}\limits_{#2}^{#3}}}
\def\binom#1#2{{{#1}\choose{#2}}}

\def\C{{\mathbb{C}}}
\def\D{{\mathbb{D}}}

\def\E{{\mathbb{E}}}

\def\End{{\rm  End}}
\def\ee{\end{equation}}

\def\epsilon{\varepsilon}

\def\gl{{\mathfrak{gl}}}

\def\id{{\rm{id}}}
\def\Id{{\rm{Id}}}

\def\lra{\longrightarrow}
\def\M{{\mathbb M}}
\def\N{{\mathbb N}}

\def\P{{\mathbb{P}}}
\def\PP{{\mathcal{P}}}
\def\pf{\noindent \textbf{Proof -- }}

\def\qed{\hfill\hbox{\vrule\vbox to 2mm{\hrule width 2mm\vfill\hrule}\vrule}  \\}
\def\R{{\mathbb{R}}}
\def\Ram{{\mathcal{R}}}
\def\S{{\mathfrak{S}}}

\def\so{{\mathfrak{so}}}
\def\sp{{\mathfrak{sp}}}
\def\su{{\mathfrak{su}}}
\def\St{{\rm St}}

\def\text{\mbox}
\def\tr{{\mathop{\rm Tr}}}
\def\u{{\mathfrak{u}}}
\def\U{{\mathcal{U}}}
\def\Un{{{U}(N)}}
\def\UC{{\mathbb{U}}}

\def\Var{{\rm Var }}

\def\Z{{\mathbb{Z}}}

\title{Schur-Weyl duality and the heat kernel measure on the unitary group}
\author{Thierry L\'EVY}
\date{January 2008}

\begin{document}
\maketitle

\begin{abstract}
We investigate a relation between the Brownian motion on the unitary group and the most natural
random walk on the symmetric group, based on Schur-Weyl duality. We use this relation to establish a
convergent power series expansion for the expectation of a product of traces of powers of a random
unitary matrix under the heat kernel measure. This expectation turns out to be the generating
series of certain paths in the Cayley graph of the symmetric group. Using our expansion, we
recover asymptotic results of Xu, Biane and Voiculescu. We give an interpretation of our main
expansion in terms of random ramified coverings of a disk. 
\end{abstract}

\section{Introduction} In this paper, we are concerned with the asymptotics of large random
unitary matrices distributed according to the heat kernel measure. This problem has been
studied first about ten years ago by P. Biane \cite{Biane} and F. Xu \cite{Xu}. It shares some
similarities with the case of unitary matrices distributed under the Haar measure, studied by B.
Collins and P. \'Sniady \cite{Collins,CollinsSniady}. The origin of our interest in
this problem is the hypothetical existence of a large $N$ limit to the two-dimensional $\Un$ Yang-Mills
theory. This limit has been investigated by physicists, in particular by V.
Kazakov and V. Kostov \cite{KK} and by D. Gross, in collaboration with W. Taylor \cite{GrossTaylor2}, A.
Matytsin \cite{GrossMatytsin} and R. Gopakumar \cite{GopakumarGross}. In \cite{Singer}, I. Singer has given the name of "Master field" to
this limit, which still has to be constructed. A. Sengupta has described in \cite{Sengupta} the
relationship between Yang-Mills theory and large unitary matrices. We refer
the interested reader to this paper and will not develop this motivation further.
Sengupta's work also contains some results whose study was at the origin of this paper (see
Proposition \ref{prop formula U} and the discussion thereafter). 

Our approach relies on the fact that the Schur-Weyl duality determines a (non-bijective)
correspondence
between conjugation-invariant objects on the unitary group on one hand and on the symmetric group
on the other hand. To be specific, let $n,N\geq 1$ be integers. Let $\rho_{n,N} :\S_{n}\times U(N)
\lra GL((\C^N)^{\otimes n})$ be the classical representation. The set $\PP_{n,N}$ of partitions of
$n$ with at most $N$ parts indexes irreducible representations of both $\S_{n}$ and $U(N)$. If
$\lambda$ is such a partition, let $\chi^\lambda$ (resp. $\chi_{\lambda}$) denote the
corresponding character on $\S_{n}$ (resp. $U(N)$). Let $Z$ be an element of the centre of
$\C[\S_{n}]$. Let $D$ be a conjugation-invariant distribution on $U(N)$. Then the equalities 
\begin{equation}\label{corresp SW}
\forall \lambda \in \PP_{n,N} , \;
\frac{\chi^{\lambda}(Z)}{\chi^{\lambda}(\id)}=\frac{\chi_{\lambda}(D)}{\chi_{\lambda}(I_{N})},
\end{equation}
where $\chi_{\lambda}(D)=D\chi_{\lambda}$, imply $\rho_{n,N}(Z\otimes 1)=\rho_{n,N}(1\otimes
D)$. The main observation, implicit in \cite{GrossTaylor2}, is the following: the element
$Z=-\frac{Nn}{2}-\sum_{1\leq k<l\leq n} (k l) \in \C[\S_{n}]$ and the distribution $D$
on $U(N)$ defined by $D\varphi=\frac{1}{2}\Delta_{U(N)} \varphi(I_{N})$, where $\Delta_{U(N)}$ is the
Laplace operator, satisfy (\ref{corresp SW}). Now $Z$ is, up to an additive constant, the generator
of the most natural random walk on $\S_{n}$ and it follows from this discussion that this random
walk is closely related to the Brownian motion on the unitary group. This relation is stated
precisely and proved in Section 2. It is also partially generalized to the orthogonal and
symplectic groups.

In Section 3, we prove our main result, which is the following.

\begin{theorem}[see also Thm \ref{main expansion}]\label{thm intro}
Let $N, n\geq 1$ be integers. Let $(B_{t})_{t\geq 0}$ be a Brownian motion on $\Un$ starting at the
identity and corresponding to the scalar product $(X,Y)\mapsto -\tr(XY)$ on
$\u(N)$. Let $\sigma$ be an element of $\S_{n}$.  Let $m_{1},\ldots,m_{r}$ denote the lengths of the
cycles of $\sigma$. Then, for all $t\geq 0$, we have the following series expansion:
\begin{equation}\label{main eqn intro}
\E[\tr_{N}(B_{\frac{t}{N}}^{m_{1}})\ldots \tr_{N}(B_{\frac{t}{N}}^{m_{r}})] =  e^{-\frac{nt}{2}}
\sum_{k,d=0}^{+\infty}\frac{(-1)^k t^k}{k!N^{2d}}S(\sigma,k,d).
\end{equation}
For all $T\geq 0$, this expansion converges uniformly on $(N,t)\in \N^* \times [0, T]$.
\end{theorem}

The coefficients $S(\sigma,k,d)$ count paths in the Cayley graph of the
symmetric group $\S_{n}$. More specifically, we consider the Cayley graph of $\S_{n}$ generated by
all transpositions. For all $\pi\in \S_{n}$, we denote by $|\pi|$ the graph distance between
$|\pi|$ and the identity. Then $S(\sigma,k,d)$ is the number of paths starting at $\sigma$
of length $k$ and finishing at a point $\pi$ such that $|\pi|=|\sigma|-(k-2d)$. In particular,
$S(\sigma,k,d)=0$ if $|k-2d|\geq n$ : for each $d\geq 0$, the contribution of order $N^{-2d}$ is
a polynomial in $t$.

The coefficients $S(\sigma,k,d)$ depend on $\sigma$ only through its conjugacy class and can be
expressed in terms of the representations of the symmetric group. In fact, Theorem \ref{thm intro}
can be proved directly using the representation theory of the unitary and symmetric groups. We
present this proof in Section 4. It is more systematical than the proof presented in Section 3 and
should be easier to generalize, as also suggested by the work of Gross and Taylor \cite{GrossTaylor2}.

The tools of representation theory allow us, in Section 5, to compute $S(\sigma,k,d)$ when
$\sigma$ is a cycle of length $n$. The expression involves Stirling numbers and it could hardly be
called simple. Nevertheless, it allows us to count for all integer $p$ the number of
ways to write the cycle $(1\ldots n)$ in $\S_{n}$ as a product of $p$ transpositions.

In Section 6, we use our expansion to describe the asymptotic distribution of unitary matrices
under the heat kernel measure as their size tends to infinity, thus recovering a result of P. Biane
\cite{Biane}. We also recover a
result of F. Xu \cite{Xu} on the asymptotic factorization of the expected values of products of
traces. In order to describe the asymptotic distribution, we must compute the coefficients
$S(\sigma,k,0)$. The factorization result mentioned above reduces the problem to the case where
$\sigma$ is an $n$-cycle. Unfortunately, the expression of $S((1\ldots n),k,0)$ obtained in Section 5
is not obviously equal to what it should be according to Biane's results. Thus, we compute this
coefficient in a different way by using the relations between the geometry of the Cayley graph of
the symmetric group and the lattice of non-crossing partitions. Then, in Section 7, we apply the
same ideas related to non-crossing partitions and use Speicher's criterion of freeness to prove the
asymptotic freeness of independent unitary matrices under the heat kernel measure. 

Finally, in Section 8, we give an interpretation of our formula in terms of random ramified
coverings over a disk, thus proving a formula described by Gross and Taylor \cite{GrossTaylor2}. We
define a probability measure on a certain set of ramified coverings over the disk and prove that
the expectation computed in Theorem \ref{thm intro} is the integral of a
simple function - essentially $N$ raised to a power equal to the Euler characteristic of the
total space of the covering - against this measure. From this point of view, our expansion deserves to be called a genus
expansion. 

It is a pleasure to thank Philippe Biane for several enlightening conversations.

\section{Probabilistic aspects of Schur-Weyl duality}

In this first section, we establish formulae which relate the heat kernel measures on $\Un$, $S\Un$,
$SO(N)$ and $Sp(N)$, to natural random walks in the symmetric group and the
Brauer monoid. 

\subsection{The unitary group}\label{unitary}

Let $n$ and $N$ be two positive integers. There is a natural action of each of the groups $\Un$
and $\S_{n}$ on the vector space $(\C^{N})^{\otimes n}$, defined as follows: for all $U\in \Un$,
$\sigma\in \S_{n}$ and $x_{1},\ldots,x_{n}\in \C^N$, we set
\begin{eqnarray}
U \cdot (x_{1}\otimes \ldots \otimes x_{n})&=& Ux_{1}\otimes \ldots \otimes U x_{n}, \nonumber \\
\sigma \cdot (x_{1}\otimes \ldots \otimes x_{n})&=& x_{\sigma^{-1}(1)}\otimes \ldots
\otimes x_{\sigma^{-1}(n)}. \label{def action}
\end{eqnarray}
It is a basic observation that these actions commute to each other. In particular, they determine an
action $\rho_{n,N}$ of $\S_{n}\times \Un$ on $(\C^{N})^{\otimes n}$ by
$$\rho_{n,N}(\sigma,U) (x_{1}\otimes \ldots \otimes x_{n}) = Ux_{\sigma^{-1}(1)}\otimes \ldots
\otimes U x_{\sigma^{-1}(n)}.$$

\begin{definition} \label{power sums 1}
Let $M_{1},\ldots,M_{n}$ be $N\times N$ complex matrices. Let $\sigma$ be an
element of $\S_{n}$. We denote by $p_{\sigma}^{st}(M_{1},\ldots,M_{n})$ the complex number
\begin{eqnarray*}
p_{\sigma}^{st}(M_{1},\ldots,M_{n})&=&\tr_{(\C^{N})^{\otimes n}}\left((M_{1}\otimes \ldots \otimes
M_{n}) \circ \rho_{n,N}(\sigma,I_N)\right)\\
&=&\prod_{\substack{c=(i_{1}\ldots i_{r})\\ {\text{\rm \scriptsize cycle of }}\sigma}}\tr(M_{i_{1}}\ldots
M_{i_{r}}).
\end{eqnarray*}
We set $p_{\sigma}^{st}(M)=p^{st}_{\sigma}(M,\ldots,M)$.
\end{definition}

The upper index {\it st} indicates that we use the standard trace rather than the normalized one
in the definition. The letter $p$ stand for "power sums", since $p^{st}_{\sigma}(M)$, as a
symmetric function of the eigenvalues of $M$, is the product of power sums corresponding to the
partition determined by $\sigma$. Observe that, by definition, the character of the representation
$\rho_{n,N}$ is the function $\chi_{\rho_{n,N}}(\sigma,U)= p^{st}_{\sigma}(U)$.

The core result of Schur-Weyl duality is that the two subalgebras of $\End((\C^{N})^{\otimes n})$
generated respectively by the actions of $\Un$ and $\S_{n}$ are each other's commutant. Let us explain
why this makes a relation between the Brownian motion on $U(N)$ and some element of the
centre of the group algebra of $\S_{n}$ unavoidable.

Let $\u(N)$ denote the Lie algebra of $\Un$, which consists of the $N\times N$ anti-Hermitian
complex matrices. Let  $\U(\u(N))$ denote the enveloping algebra of $\u(N)$, which is canonically
isomorphic to the algebra of left-invariant differential operators on $\Un$. Let also $\C[\S_{n}]$
denote the group algebra of $\S_{n}$. The representation $\rho_{n,N}$ determines a homomorphism of
associative algebras  $\C[\S_{n}]\otimes \U(\u(N))\lra \End((\C^{N})^{\otimes n})$. 
The centre $\mathcal{Z}(\u(N))$ of $\U(\u(N))$ is the space of bi-invariant differential
operators on $\Un$. Since $\rho_{n,N}(1\otimes \mathcal{Z}(\u(N)))$ commutes with
$\rho_{n,N}(1,U)$ for every $U\in \Un$, the Schur-Weyl duality asserts in particular that
$$\rho_{n,N}(1\otimes \mathcal{Z}(\u(N))) \subset \rho_{n,N}(\C[\S_{n}]\otimes 1).$$
We are primarily interested in the Laplace operator, which is defined as follows. The
$\R$-bilinear form $\langle X,Y\rangle=\tr(X^*Y)=-\tr(XY)$ is a scalar product on $\u(N)$. Let $(X_{1},\ldots,X_{N^2})$
be an orthonormal basis of $\u(N)$. Identifying the elements of $\u(N)$ with left-invariant vector
fields on $\Un$, thus with first-order differential operators on $U(N)$, the Laplace operator $\Delta_{\Un}$
is the differential operator
$\sum_{i=1}^{N^2} X_{i}^{2}$. It corresponds to the Casimir element $\sum_{i=1}^{N^2} X_{i}\otimes
X_{i}$ of the enveloping algebra of $\u(N)$. This element is central and does not depend on the
choice of the orthonormal basis. Hence, $\Delta_{\Un}$ is well defined and bi-invariant. The
discussion above shows that, in the representation $\rho_{n,N}$, the
Laplace operator of $\Un$ can be expressed as an element of $\C[\S_{n}]$. This is exactly what the
main formula of this section does, in an explicit way.

Let $T_{n}$ be the subset of $\S_{n}$ consisting of all transpositions. We set $\Delta_{\S_{n}}=-\frac{n(n-1)}{2}+ \sum_{\tau
\in T_{n}} \tau$. The formula for the unitary group is the following.

\begin{proposition} \label{prop formula U} For all integers, $n,N\geq 1$, one has
\begin{equation}\label{formula U}
\rho_{n,N}\left(\Delta_{\S_{n}}\otimes 1 + 1\otimes \frac{1}{2}\Delta_{\Un}\right)= - \frac{Nn+n(n-1)}{2}.
\end{equation}
\end{proposition}

Before we prove this formula, let us derive some of its consequences.

\begin{proposition} \label{lap newt st} For each $\sigma\in \S_{n}$, the function
$p^{st}_{\sigma}:\Un\lra \C$ satisfies the following relation:
\begin{equation}
\frac{1}{2}\Delta_{\Un}p^{st}_{\sigma} = -\frac{Nn}{2}p^{st}_{\sigma} - \sum_{\tau
\in T_{n}}p^{st}_{\sigma\tau}. \label{delta p st}
\end{equation}
More generally, let $M_{1},\ldots,M_{n}$ be arbitrary $N\times N$ matrices. Then, regarding
$p^{st}_{\sigma}(M_{1}U,\ldots,M_{n}U)$ as a function of $U\in \Un$, one has
\begin{equation}
\hskip -0.3cm \frac{1}{2}\Delta_{\Un}p^{st}_{\sigma}(M_{1}U,\ldots,M_{n}U)=-\frac{Nn}{2} 
p^{st}_{\sigma}(M_{1}U,\ldots,M_{n}U) - \sum_{\tau\in
T_{n}}p^{st}_{\sigma\tau}(M_{1}U,\ldots,M_{n}U). \label{delta p st M}
\end{equation}
\end{proposition}

\pf Recall that $p^{st}_{\sigma}(M_{1}U,\ldots,M_{n}U)=\tr\left((M_{1}\otimes\ldots
\otimes M_{n})\circ \rho_{n,N}(\sigma,U)\right)$. Let us use the shorthand notation $M=M_{1}\otimes
\ldots \otimes M_{n}$. We have
\begin{eqnarray*}
\frac{1}{2}\Delta_{\Un}p^{st}_{\sigma}(M_{1}U,\ldots,M_{n}U)&=&\tr(M\circ
\rho_{n,N}(\sigma,U)\circ \rho_{n,N}(1\otimes \frac{1}{2}\Delta_{\Un}))\\
&&\hskip -3.5cm = -\frac{Nn+n(n-1)}{2}p^{st}_{\sigma}(M_{1}U,\ldots,M_{n}U)  - \tr(M\circ \rho_{n,N}(\sigma,U)\circ \rho_{n,N}(\Delta_{\S_{n}}\otimes 1)).
\end{eqnarray*}
The result follows immediately from the definition of $\Delta_{\S_{n}}$. \qed

The function $p^{st}_{\sigma}$ depends only on the cycle structure of $\sigma$. In concrete terms,
if the lengths of the cycles of $\sigma$ are $m_{1},\ldots,m_{r}$, then
$p^{st}_{\sigma}(U)=\tr(U^{m_{1}})\ldots \tr(U^{m_{r}})$. This redundant labelling is however
nicely adapted to our problem, as equation (\ref{delta p st}) shows. Let us spell out the right
hand side of this equality.  The permutation $\sigma$ being fixed, the cycle structure of
$\sigma\tau$ depends on the two points exchanged by the transposition $\tau$. If they belong to the
same cycle of $\sigma$, then this cycle is split into two cycles. A cycle of length $m$ can be
split into a cycle of length $s$ and a cycle of length $m-s$ by $m$ distinct transpositions,
unless $m=2s$, in which case only $\frac{m}{2}$ of these transpositions are distinct. If on the
contrary the points exchanged by $\tau$ belong to two distinct cycles of $\sigma$, these two cycles
are merged into a
single cycle. Two cycles of lengths $m$ and $m'$ can be merged by $mm'$ distinct permutations.
Altogether, we find the following equation, which was already present in papers of Xu \cite{Xu}
and Sengupta \cite{Sengupta}.

\begin{eqnarray*}
\Delta_{\Un}\left(\tr(U^{m_{1}})\ldots\tr(U^{m_{r}})\right)&=&
-Nn \tr(U^{m_{1}})\ldots\tr(U^{m_{r}}) \\
&& \hskip -2.5cm + \sum_{i=1}^{r}  m_{i} \tr(U^{m_{1}})\ldots \widehat{\tr(U^{m_{i}})}\ldots
\tr(U^{m_{r}})\;  \sum_{s=1}^{m_{i}-1} \tr(U^{s})\tr(U^{m_{i}-s})\\
&&\hskip -2.5cm +\sum_{\substack{i,j=1, i\neq j}}^{r} m_{i}m_{j} \tr(U^{m_{1}})\ldots \widehat{\tr(U^{m_{i}})}\ldots
\widehat{\tr(U^{m_{j}})}\ldots \tr(U^{m_{r}})\; \tr(U^{m_{i}+m_{j}}).
\end{eqnarray*}

A remarkable feature of (\ref{formula U}) is the fact that
the element of $\C[\S_{n}]$ which appears has coefficients of the same sign on
the elements which are not the identity. Hence, up to an additive constant, it can be interpreted as
the generator of a Markov chain on $\S_{n}$. This leads us to the following simple probabilistic
interpretation of (\ref{formula U}). 

Let us introduce the standard random walk on the Cayley graph of the symmetric group
generated by the set of transpositions. It is the continuous-time Markov chain on $\S_{n}$ with
generator $\Delta_{\S_{n}}$, that is, the chain which jumps at rate  ${n}\choose{2}$ from its current
position $\sigma$ to $\sigma\tau$, where $\tau$ is
chosen uniformly at random among the ${n}\choose{2}$ transpositions of $\S_{n}$. 

If $\sigma$ is a permutation, we denote by $\ell(\sigma)$ the number of cycles of $\sigma$. For
example, $\tau$ is a transposition if and only if $\ell(\tau)=n-1$.

\begin{proposition} \label{BM SRW}Let $N, n\geq 1$ be integers. Let $(B_{t})_{t\geq 0}$ be a Brownian
motion on $\Un$ starting at the identity and corresponding to the scalar product $(X,Y)\mapsto -\tr(XY)$ on
$\u(N)$. Let $(\pi_{t})_{t\geq 0}$ be a standard random walk on the
Cayley graph of the symmetric group $\S_{n}$, independent of $(B_{t})_{t\geq 0}$. Then the process $\left(e^{\frac{Nn+n(n-1)}{2}t}p^{st}_{\pi_t}(B_t)\right)_{t\geq 0}$ is a martingale. In particular,
\begin{equation}\label{prob U}
\E\left[p^{st}_{\pi_{t}}(B_{t})\right]=e^{-\frac{Nn+n(n-1)}{2}t}\;  \E\left[N^{\ell(\pi_{0})}\right].
\end{equation}
More generally,  let $M_{1},\ldots,M_{n}$ be arbitrary $N\times N$ complex matrices. Then the stochastic process $\left(e^{\frac{Nn+n(n-1)}{2}t}p^{st}_{\pi_t}(M_1B_t,\ldots,M_n B_t)\right)_{t\geq 0}$ is a martingale and
\begin{equation}\label{prob M U}
\E\left[p^{st}_{\pi_{t}}(M_{1}B_{t},\ldots,M_{n}B_{t})\right]=e^{-\frac{Nn+n(n-1)}{2}t}\; 
\E\left[p^{st}_{\pi_{0}}(M_{1},\ldots,M_{n})\right].
\end{equation}
\end{proposition}

\pf The process $(\pi,B)$ is a Markov process on $\S_{n}\times \Un$ with generator
$\Delta_{\S_{n}}\otimes 1+1\otimes \frac{1}{2}\Delta_{\Un}$. Consider the function $p:\S_{n}\times \Un \to \C$ defined by $p(\sigma,U)=p^{st}_\sigma(M_1 U,\ldots,M_n U)$. By Proposition \ref{lap newt st}, this function satisfies the relation
$$\left(\Delta_{\S_{n}}\otimes 1+1\otimes \frac{1}{2}\Delta_{\Un}\right) p = -\frac{Nn+n(n-1)}{2} p.$$
The fact that $\left(e^{\frac{Nn+n(n-1)}{2}t}p^{st}_{\pi_t}(M_1B_t,\ldots,M_n B_t)\right)_{t\geq 0}$ is a martingale follows immediately. The last assertion follows from the fact that $B_{0}=I_{N}$ a.s. \qed

Let us turn to the proof of Proposition \ref{prop formula U}. \\

{\noindent \textbf{Proof of Proposition \ref{prop formula U} -- }} The action of $\u(N)$ on
$(\C^{N})^{\otimes n}$ extends by complexification to $\gl(N,\C)=\u(N)\oplus i\u(N)$. Let
$(X_{1},\ldots,X_{N^2})$ be a real basis of $\u(N)$. It is also a complex basis of $\gl(N,\C)$. 
Define a $N\times N$ matrix $g$ by  $g_{ij}=-\tr(X_{i}X_{j})$. Since $-\tr(\cdot\;\cdot)$ is non-degenerate on $\gl_{N}(\C)$, the 
matrix $g$ has an inverse $g^{-1}$, the entries of which we denote by $g^{ij}$. Then it is easy to
check that the element  $\sum_{i,j=1}^{N^2} g^{ij} X_{i}\otimes X_{j}$ of the enveloping algebra
is independent of the choice of the basis. Of course, by choosing our original basis of $\u(N)$
orthonormal, we find that this element is simply $\Delta_{\Un}$. 

In order to compute $\rho_{n,N}(1\otimes \Delta_{\Un})$, we prefer to use 
another complex basis of $\gl(N,\C)=\M_{N}(\C)$, namely the canonical basis $(E_{ij})_{i,j\in
\{1,\ldots,N\}}$. For this basis,  $g_{ij,kl}=-\delta_{jk}\delta_{il}$ and $g=g^{-1}$. Hence, in the enveloping algebra of $\gl(N,\C)$, $\Delta_{U(N)}= -\sum_{i,j=1}^{N} E_{ij}\otimes E_{ji}$.

First, notice that $\rho_{n,N}(1\otimes E_{ij})(x_{1}\otimes \ldots \otimes x_{n})=\sum_{k=1}^{n} x_{1}\otimes
\ldots \otimes E_{ij}(x_{k})\otimes\ldots \otimes x_{n}.$ Hence,
\begin{eqnarray*}
\hskip -0cm \rho_{n,N}\left(1\otimes \sum_{i,j=1}^{N}E_{ij}\otimes E_{ji}\right) &=& 2\sum_{i,j=1}^{N}\sum_{1\leq
k<l\leq n}
\Id^{\otimes
k-1}\otimes E_{ij}\otimes \Id^{\otimes l-k-1} \otimes E_{ji} \otimes \Id^{\otimes n-l-1} +\\
&& + \sum_{i,j=1}^{N} \sum_{k=1}^{n} \Id^{\otimes k-1}\otimes E_{ii} \otimes \Id^{\otimes n-k-1}.
\end{eqnarray*}
The last term is simply $Nn$ times the identity. For the first part of the right hand side,
observe that $\sum_{i,j=1}^{N} E_{ij}\otimes E_{ji}\in \End((\C^{N})^{\otimes 2})$ is the
transposition operator $x\otimes y\mapsto y\otimes x$, that is, the operator $\rho_{2,N}((12),I_{N})$.
Finally, we have found that
$$-\rho_{n,N}(1\otimes \Delta_{\Un})= Nn\; \Id + \sum_{1\leq k\neq l \leq n} \rho_{n,N}((kl),I_{N}).$$
The result follows.\qed

The results of this section still hold, after a minor modification, when $\Un$ is replaced by
$S\Un$. Indeed, the orthogonal complement of $\su(N)$ in $\u(N)$ is the line generated by
$\frac{i}{\sqrt{N}}I_{N}$. Since $\rho_{n,N}(1\otimes (I_{N}\otimes I_{N}))=n^2\Id$, the Casimir
operator of $\su(N)$ satisfies the relation
$$\rho_{n,N}(1\otimes \Delta_{S\Un}) =\rho_{n,N}(1\otimes \Delta_{\Un})+\frac{n^2}{N}\Id.$$
This modifies only the exponential factors in (\ref{prob U}) and (\ref{prob M U}).

We will explore further consequences of Proposition \ref{BM SRW} in the rest of the
paper. For the moment, we derive similar results for the orthogonal and symplectic group.

\subsection{The orthogonal group} 

Let us consider the action of $SO(N)$ on $(\C^{N})^{\otimes n}$ defined by analogy with (\ref{def
action}). The action of $\S_{n}$ still commutes to that of $SO(N)$, but, unless $n=1$, the
subalgebra of $\End((\C^{N})^{\otimes n})$ generated by the image of $\C[\S_{n}]$ is strictly
smaller than the commutant of the image of $SO(N)$. Let us review briefly the operators which are
classically used to describe this commutant. We denote by $\{e_{1},\ldots,e_{N}\}$ the canonical
basis of $\C^{N}$.

\begin{definition} Let $\beta$ be a partition of $\{1,\ldots,2n\}$ into pairs. Define 
$\rho_{n,N}(\beta)\in \End((\C^{N})^{\otimes n})$ by setting, for all $i_{1},\ldots,i_{n}\in \{1,\ldots,N\}$,
$$\rho_{n,N}(\beta)(e_{i_{1}}\otimes \ldots \otimes e_{i_{n}})= \sum_{i_{n+1},\ldots,i_{2n}\in \{1,\ldots,N\}}
\prod_{\{k,l\}\in \beta} \delta_{i_{k}i_{l}}\; e_{i_{n+1}}\otimes \ldots \otimes e_{i_{2n}}.$$
\end{definition}

Observe that the partition $\{\{1,n+1\},\ldots,\{n,2n\}\}$ is sent to the identity operator by
$\rho_{n,N}$. Let $B_{n}$ denote the set of partitions of $\{1,\ldots,2n\}$ into pairs. The
composition of
the operators $\rho_{n,N}(\beta)$ corresponds to a monoid structure on $B_{n}$ which is easiest to
understand on a picture. An element of $B_{n}$ is represented in a box with $n$ dots on its top
edge and $n$ dots on its bottom edge. The dots on the top are labelled from $1$ to $n$, from the
left to the right. The dots on the bottom are labelled from $n+1$ to $2n$, from the left to the
right too. A pairing is then simply represented by $n$ chords which join the appropriate dots.
Multiplication of pairings is done in the intuitive topological way by superposing boxes and, if
necessary, removing the closed loops which have appeared.

\begin{figure}[h!]
\begin{center}
\scalebox{0.8}{\includegraphics{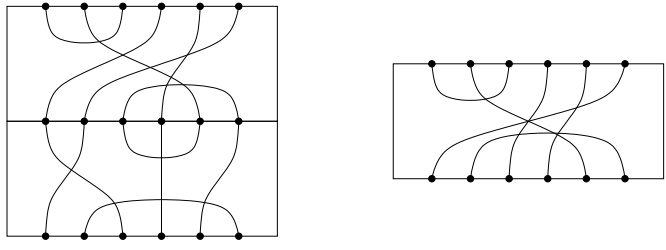}}
\caption{Multiplication of two diagrams in the Brauer monoid.}
\end{center}
\end{figure}

The monoid $B_{n}$ is called the Brauer monoid and its elements are called Brauer
diagrams. The group $\S_{n}$ is naturally a submonoid\footnote{In fact, $\S_{n}\subset B_{n}$ is
exactly the subset of invertible elements.
Indeed, for $\beta\in B_{n}$, let $T(\beta)$ be the set of pairs $\{k,l\}\in \beta$ such that
$1\leq k,l\leq n$. In words, $T(\beta)$ is the set of chords in the diagram of $\beta$ which
join two dots on the top edge of the box. It is clear that $T(\beta_{1}\beta_{2})\supset
T(\beta_{1})$ for all $\beta_{1},\beta_{2}\in B_{n}$. Hence, $T(\beta)$ must be empty for $\beta$
to be invertible. More generally, it is not difficult to check that, given $\beta$ and $\beta'$ in
$B_{n}$, there exists $\beta''\in B_{n}$ such that $\beta\beta''=\beta'$ if and only if 
$T(\beta)\subset T(\beta')$.} of $B_{n}$, by the identification of a
permutation $\sigma$ with the pairing $\{\{1,\sigma(1)+n\},\ldots,\{n,\sigma(n)+n\}\}$. The
identification of $\S_{n}$ with a subset of $B_{n}$ is compatible with our previous definition
of $\rho_{n,N}$ in the sense that $\rho_{n,N}(\sigma)$ is the same if we consider
$\sigma$ as a permutation or as a Brauer diagram.

The correct statement of Schur-Weyl duality in the present context is that the subalgebras of $\End((\C^{N})^{\otimes
n})$ generated by $SO(N)$ and $B_{n}$ are each other's commutant (see \cite{GoodmanWallach}). Let $\rho_{n,N}$ denote the
morphism of monoids
$$\rho_{n,N}:B_{n}\times SO(N)\lra GL((\C^{N})^{\otimes n}).$$
Just as in the unitary case, this action determines a morphism of associative algebras
$\rho_{n,N}:\C[B_{n}]\otimes \U(\so(N))\lra \End((\C^{N})^{\otimes n})$.

By analogy to the unitary case, let us define "power sums" functions associated to Brauer diagrams.
Given $\beta\in B_{n}$ and $M_{1},\ldots,M_{n}\in \M_{N}(\C)$, set
$$p^{st}_{\beta}(M_{1},\ldots,M_{n})=\tr((M_{1}\otimes \ldots \otimes M_{n}) \circ
\rho_{n,N}(\beta)).$$
In particular, the character of $\rho_{n,N}$ is given by $\chi_{\rho_{n,N}}(\beta,R)=p^{st}_{\beta}(R)$. 

The number $p^{st}_{\beta}(M_{1},\ldots,M_{n})$ is a product of traces of words in the matrices
$M_{1},{}^{t}M_{1}$ $,\ldots,$ $M_{n},{}^{t}M_{n}$. Let us describe in more detail how to
compute
$p^{st}_{\beta}(I_{N})$. Let $\beta$ be a Brauer diagram. Consider the graph with vertices $\{1,\ldots,n\}$ and unoriented edges
$\{k,l\}$, where $k$ and $l$ are such that there exist $k'\in\{k,k+n\}$ and $l'\in\{l,l+n\}$ with
$\{k',l'\}\in\beta$. This is the graph obtained by identifying the top edge with the bottom edge
in the graphical representation of $\beta$. Then each vertex has degree $2$ in this graph. Hence,
it is a union of disjoint unoriented cycles. If $\beta$ belongs to $\S_{n}\subset B_{n}$, this
cycle structure is of course that of $\beta$ as a permutation, apart from the orientation which is
lost. In general, let $\ell(\beta)$ denote the number of cycles in this graph. Then
$p^{st}_{\beta}(I_{N})=N^{\ell(\beta)}$.

Let us define an element of $\C[B_{n}]$ as follows. Given $k$ and $l$ two integers such that
$1\leq k<l\leq n$, we define the element $\langle kl\rangle$ of $B_{n}$ as the following pairing:
$$\langle kl\rangle=\{\{k,l\},\{n+k,n+l\}\}\cup \bigcup_{i\in\{1,\ldots,n\}-\{k,l\}} \{\{i,n+i\}\}.$$

\begin{figure}[h!]
\begin{center}
\scalebox{0.8}{\includegraphics{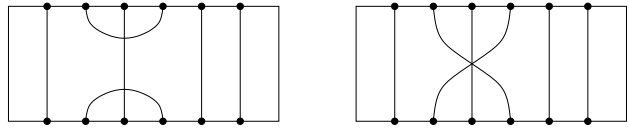}}
\caption{The elements $\langle 24 \rangle$ and $(24)$ of $B_{6}$.}
\end{center}
\end{figure}

Let $C_{n}$ be the subset of $B_{n}$ consisting of all the element of the form $\langle kl
\rangle$. We now define $\displaystyle \Delta_{B_{n}}=-\frac{n(n-1)}{2}+\sum_{\alpha\in C_{n}}
\alpha$.
Thanks to the inclusion $\S_{n}\subset B_{n}$, we still see $\Delta_{\S_{n}}$ as an element of
$\C[B_{n}]$. The formula for the orthogonal group is the following.

\begin{proposition} \label{prop formula O} For all integers, $n,N\geq 1$, one has
\begin{equation}\label{formula O}
\rho_{n,N}\left(\Delta_{\S_{n}}\otimes 1 + 1\otimes \Delta_{SO(N)}\right)= -
\frac{(N-1)n}{2}+\rho_{n,N}\left(\Delta_{B_{n}}\otimes 1\right).
\end{equation}
\end{proposition}

\pf The computation is very similar to that we made in the unitary case. Endow $\so(N)$ with the
scalar product $\langle X, Y\rangle=-\tr(XY)$. The basis $(A_{ij})_{1\leq i<j \leq N}$, with $A_{ij}=E_{ij}-E_{ji}$,
is orthogonal and $\langle A_{ij},A_{ij}\rangle=2$ for all $i<j$. Hence, $\Delta_{SO(N)}=\frac{1}{2}\sum_{1\leq i<j\leq N}
A_{ij}\otimes A_{ij}$. We have
\begin{eqnarray*}
\hskip -0.7cm\rho_{n,N}(1\otimes \Delta_{SO(N)})&=&\sum_{1\leq k<l\leq n}\sum_{1\leq i<j\leq
N}\Id^{\otimes k-1} \otimes A_{ij}\otimes \Id^{\otimes l-k-1}\otimes
A_{ij}\otimes \Id^{\otimes n-l} +\\
&& \hskip 0cm + \frac{1}{2}\sum_{k=1}^{n}\sum_{1\leq i<j\leq N} \Id^{\otimes k-1}\otimes A_{ij}^{2}\otimes
\Id^{\otimes n-k}\\
&=& \sum_{1\leq k<l\leq n}\sum_{i,j=1}^{N}\Id^{\otimes k-1} \otimes E_{ij}\otimes \Id^{\otimes l-k-1}
\otimes E_{ij} \otimes \Id^{\otimes n-l}\\
&&\hskip -0cm - \sum_{1\leq k<l\leq n} \sum_{i,j=1}^{N}\Id^{\otimes k-1} \otimes E_{ij}\otimes
\Id^{l-k-1} \otimes E_{ji}\otimes \Id^{n-l} -\frac{(N-1)n}{2} \Id\\
&=&\sum_{1\leq k<l\leq n} \rho_{n,N}((\langle kl\rangle-(kl))\otimes 1) -\frac{(N-1)n}{2} \Id.
\end{eqnarray*}
The result follows.\qed

The following proposition is proved just as Proposition \ref{lap newt st}.

 \begin{proposition} \label{lap newt st O} For all $\beta\in B_{n}$, the following relation holds:
\begin{equation}
\Delta_{SO(N)}p^{st}_{\beta} = -\frac{(N-1)n}{2}p^{st}_{\beta} - \sum_{\tau
\in T_{n}}p^{st}_{\beta\tau} + \sum_{\alpha \in C_{n}}p^{st}_{\beta\alpha}. \label{delta p st O}
\end{equation}
More generally, let $M_{1},\ldots,M_{n}$ be arbitrary $N\times N$ matrices. Then, regarding
$p^{st}_{\beta}(M_{1}R,\ldots,M_{n}R)$ as a function of $R\in SO(N)$,
\begin{eqnarray}
 \Delta_{SO(N)}p^{st}_{\beta}(M_{1}R,\ldots,M_{n}R)&=&-\frac{(N-1)n}{2} 
p^{st}_{\beta}(M_{1}R,\ldots,M_{n}R) \nonumber \\
&&\hskip 0cm -\sum_{\tau\in T_{n}}p^{st}_{\beta\tau}(M_{1}R,\ldots,M_{n}R) +\sum_{\alpha\in C_{n}}p^{st}_{\beta\alpha}(M_{1}R,\ldots,M_{n}R). \label{delta p st M O}
\end{eqnarray}
\end{proposition}

It seems more difficult to find a probabilistic interpretation of (\ref{formula O}) than in the
unitary case, because the element of $\C[B_{n}]$ which appears does not have coefficients of the
same sign on all elements not equal to 1.

\subsection{The symplectic group} Nothing really new is needed to treat the case of the symplectic
group. Let us describe briefly the results. 

Let $J\in \M_{2N}(\C)$ denote the matrix $\left(\begin{array}{cc}0 & I_{N}\cr -I_{N} &0
\end{array}\right)$.
The symplectic group is defined by $Sp(N)=\{S\in U(2N):{}^{t}SJS=J\}$. It acts naturally on
$((\C^{2N})^{\otimes n})$. The action of the Brauer monoid needs to be slightly modified to fit
the symplectic case. If $\beta$ belongs to $B_{n}$, we define the operator $\rho_{n,2N}(\beta)$ by
setting, for all $i_{1},\ldots,i_{n}\in \{1,\ldots,2N\}$,
$$\rho_{n,2N}(\beta)(e_{i_{1}}\otimes \ldots \otimes e_{i_{n}})= \sum_{i_{n+1},\ldots,i_{2n}\in \{1,\ldots,n\}}
\prod_{\{k,l\}\in \beta} J_{i_{k}i_{l}}\; e_{i_{n+1}}\otimes \ldots \otimes e_{i_{2n}}.$$
Then we have an action $\rho_{n,2N}: B_{n}\times Sp(N)\lra \End((\C^{2N})^{\otimes n})$ and the
images of $B_{n}$ and $Sp(N)$ generate two algebras which are each other's commutant. 

The Lie algebra $\sp(N)$ is endowed with the scalar product $\langle X,Y\rangle=-\tr(XY)$ and we
denote by $\Delta_{Sp(N)}$ the corresponding Laplace operator. The main formula is the following.

\begin{proposition}\label{prop formula S}For all integers, $n,N\geq 1$, one has
\begin{equation}\label{formula 2}
\rho_{n,2N}\left(\Delta_{\S_{n}}\otimes 1 + 1\otimes 2\Delta_{Sp(N)}\right)= - (2N+1)n +\rho_{n,2N}\left(\Delta_{B_{n}}\otimes 1\right).
\end{equation}
\end{proposition}

\pf Just as in the unitary case, it is more convenient to use complexi\-fi\-cation. The Lie algebra
$\sp(N,\C)=\sp(N)\oplus i\sp(N)$ is the Lie subalgebra of $\gl(2N,\C)$ defined by the relation
${}^{t}XJ=-JX$. It consists of the matrices $\left(\begin{array}{cc}A & B\cr C &-{}^{t}A
\end{array}\right)$, where $A$ is an arbitrary $N\times N$ matrix and $B,C$ are two symmetric
$N\times N$ matrices. We use the following basis of $\sp(N,\C)$:
$$\begin{array}{rcll}
 A_{ij}& =&E_{ij}-E_{j+N ,i+N} & 1\leq  i,j\leq N\cr
 B_{ij}& =&E_{i, j+N}+E_{j, i+N} & 1\leq i<j\leq N\cr
 C_{ij}& =&E_{i+N ,j}+E_{j+N ,i}  & 1\leq i<j \leq N\cr
 D_{i}& =&E_{i, i+N} & 1\leq i \leq N\cr
 D_{i+N}&=&E_{i+N ,i} & 1\leq i \leq N.
\end{array}
$$
The bilinear form $\langle \cdot ,\cdot \rangle$ takes the following values on this basis:
\begin{eqnarray*}
\langle A_{ij},A_{ji}\rangle=-2 &\hspace{4mm}& 1\leq i,j\leq N \\
\langle B_{ij}, C_{ij}\rangle=-2 &\hspace{4mm}&1\leq i<j\leq N\\
\langle D_{i},D_{i+N}\rangle =-1 &\hspace{4mm}& 1\leq i\leq N.
\end{eqnarray*}
The other values are zero. It follows that the Casimir element of $\sp(N,\C)$ is equal to
\begin{eqnarray*}
\Delta_{\sp(N,\C)}&=&-\frac{1}{2}\sum_{1\leq i,j\leq N} A_{ij}\otimes
A_{ji}-\frac{1}{2}\sum_{1\leq i<j\leq N} (B_{ij}\otimes C_{ij} + C_{ij}\otimes B_{ij})\\
&& \hskip 5cm -\sum_{1\leq i \leq N} (D_{i}\otimes D_{i+N}+D_{i+N}\otimes D_{i}).
\end{eqnarray*}
The formula follows now by a direct computation. In order to recognize operators of the form
$\langle kl \rangle$ and $(kl)$, observe that, when $n=2$ for example,
$$ \rho_{2,2N}((12),I_{N}) = \sum_{i,j=1}^{2N} E_{ij}\otimes E_{ji},$$
\begin{eqnarray*}
\rho_{2,2N}(\langle 12 \rangle,I_{N})&=&\sum_{i,j=1}^{N} \left( E_{ij}\otimes
E_{i+N,j+N}+E_{i+N,j+N}\otimes E_{ij}\right.\\
&&\hskip 1.5cm  \left. -E_{i,j+N}\otimes E_{i+N,j}-E_{i+N,j}\otimes E_{i,j+N}\right).
\end{eqnarray*}
\qed

\begin{proposition} \label{lap newt st S} For all $\beta\in B_{n}$, the following relation holds:
\begin{equation}
2\Delta_{Sp(N)}p^{st}_{\beta} = -(2N+1)n\; p^{st}_{\beta} - \sum_{\tau
\in T_{n}}p^{st}_{\beta\tau} + \sum_{\alpha \in C_{n}}p^{st}_{\beta\alpha}. \label{delta p st S}
\end{equation}
More generally, let $M_{1},\ldots,M_{n}$ be arbitrary $2N\times 2N$ matrices. Then, regarding
$p^{st}_{\beta}(M_{1}S,\ldots,M_{n}S)$ as a function of $R\in Sp(N)$, one has
\begin{eqnarray}
2 \Delta_{Sp(N)}p^{st}_{\beta}(M_{1}S,\ldots,M_{n}S)&=&-(2N+1)n \;
p^{st}_{\beta}(M_{1}S,\ldots,M_{n}S)\nonumber \\
&& \hskip 0cm - \sum_{\tau\in T_{n}}p^{st}_{\beta\tau}(M_{1}S,\ldots,M_{n}S)  +\sum_{\alpha\in C_{n}}p^{st}_{\beta\alpha}(M_{1}S,\ldots,M_{n}S). \label{delta p st
M S}
\end{eqnarray}
\end{proposition}

\section{The power series expansion}

Let us denote by $\tr_{N}=\frac{1}{N}\tr$ the normalized trace on $\M_{N}(\C)$. Let
$M_{1},\ldots,M_{n}$ be $N\times N$ matrices. Let $\sigma$ be an element of $\S_{n}$. We denote
by $p_{\sigma}(M_{1},\ldots,M_{n})$ the number
$$p_{\sigma}(M_{1},\ldots,M_{n})=\prod_{\substack{c=(i_{1}\ldots i_{r}) \\
{\text{\rm\scriptsize cycle of }}\sigma}}\tr_{N}(M_{i_{1}}\ldots M_{i_{r}}).$$
We denote by $\ell(\sigma)$ the number of cycles of $\sigma$, so that
$p_{\sigma}=N^{-\ell(\sigma)}p^{st}_{\sigma}$.

In this section, we exploit the result of Proposition \ref{BM SRW} and derive a convergent power
series expansion of $\E\left[p_{\sigma}(B_{\frac{t}{N}})\right]$ when $B$ is a Brownian motion on $\Un$.
This expansion involves combinatorial coefficients, which count paths in the Cayley graph of
$\S_{n}$. We start by discussing these paths and introducing some notation.

\subsection{The Cayley graph of the symmetric group} Fix $n\geq 1$. The Cayley graph of $\S_{n}$
generated by $T_{n}$ can be described as follows: the vertices of this graph are the elements
of $\S_{n}$ and two permutations $\sigma_{1}$ and $\sigma_{2}$ are joined by an edge if and only if
$\sigma_{1}\sigma_{2}^{-1}$ is a transposition. It is a fundamental observation that, if $\sigma_{1}$ and $\sigma_{2}$ are joined by an edge, then $\ell(\sigma_1)$ and $\ell(\sigma_2)$ differ exactly by $1$. Indeed, multiplying a permutation by a transposition splits a cycle into two shorter cycles if the points exchanged by the transposition belong originally to the same cycle, and otherwise combines together the two cycles which contain the points exchanged by the transposition.

A finite sequence $(\sigma_{0},\ldots,\sigma_{k})$ of permutations such that $\sigma_{i}$ is joined to $\sigma_{i+1}$ by an edge for each $i\in
\{0,\ldots,k-1\}$ is called a path of length $k$. The distance between two permutations is the
smallest length of a path which joins them. This distance can be computed explicitly as follows.

Let us introduce the notation $|\sigma|=n-\ell(\sigma)$. We have $|\sigma|\in\{0,\ldots,n-1\}$ and
$|\sigma|=0$ (resp. 1, resp. $n-1$) if and only if $\sigma$ is the identity (resp. a
transposition, resp. a $n$-cycle). Other values of $|\sigma|$ do not characterize uniquely the
conjugacy class of $\sigma$. It is well-known and easy to check that $|\sigma|$ is the minimal
number of transpositions required to write $\sigma$. In other words, the graph distance
between two permutations $\sigma_{1}$ and $\sigma_{2}$ in the Cayley graph is given by $|\sigma_{1}^{-1}\sigma_{2}|$.

It turns out that the paths which play the most important role in our problem are those which tend
to get closer to the identity. Let $\gamma=(\sigma_{0},\ldots,\sigma_{k})$ be a path. Recall that, for all $i\in\{0,\ldots,k-1\}$, one has $\ell(\sigma_{i+1})=\ell(\sigma_{i})\pm 1$. We call
defect of $\gamma$ and denote by
$d(\gamma)$ the number of steps which increase the distance to the identity. In symbols,
\begin{eqnarray*}
d(\gamma)&=&\# \{i\in \{0,\ldots,k-1\} : |\sigma_{i+1}|=|\sigma_{i}|+1\}\\
&=&\#\{i\in \{0,\ldots,k-1\} :
\ell(\sigma_{i+1})=\ell(\sigma_{i})-1\}.
\end{eqnarray*}

The following lemma is straightforward.

\begin{lemma} \label{l d k}
Let $\gamma=(\sigma_{0},\ldots,\sigma_{k})$ be a path. Then $2d(\gamma)=
k-(\ell(\sigma_{k})-\ell(\sigma_{0})).$
\end{lemma}

For $\sigma,\sigma'$ in $\S_{n}$ and $k\geq 0$, let us denote by $\Pi_{k}(\sigma \to \sigma')$ the
set of paths of length $k$ which start at $\sigma$ and finish at $\sigma'$. Let us also denote 
by $\Pi_{k}(\sigma)$ the set of all paths of length $k$ starting at $\sigma$ and by
$\Pi(\sigma\to\sigma')$ the set of all paths from $\sigma$ to $\sigma'$.
Notice that the cardinality of $\Pi_{k}(\sigma)$ is equal to ${{n}\choose{2}}^{k}$. Let us
finally define the coefficients which appear in the expansion.

\begin{definition} \label{def S}
Consider $\sigma \in \S_{n}$ and two integers $k,d\geq 0$. We set
$$S(\sigma,k,d)=\#\{\gamma\in \Pi_{k}(\sigma) : d(\gamma)=d\}.$$
In words, $S(\sigma,k,d)$ is the number of paths in the Cayley graph of $\S_{n}$ starting at
$\sigma$, of length $k$ and with defect $d$.
\end{definition}

Observe that the adjoint action of $\S_{n}$ on itself determines an action of $\S_{n}$ on its
Cayley graph by automorphisms. Thus, $S(\sigma,k,d)$ depends only on the conjugacy class of
$\sigma$.

\subsection{The main expansion}

\begin{theorem} \label{main expansion} Let $N, n\geq 1$ be integers. Let $(B_{t})_{t\geq 0}$ be a
Brownian motion on $\Un$ starting at the identity and corresponding to the scalar product $(X,Y)\mapsto -\tr(XY)$ on
$\u(N)$. Let $M_{1},\ldots,M_{n}$ be arbitrary $N\times N$ complex matrices. Let $\sigma$ be an
element of $\S_{n}$. Then, for all $t\geq 0$, we have the following series expansions:
\begin{equation} 
\E\left[p_{\sigma}(M_{1}B_{\frac{t}{N}},\ldots,M_{n}B_{\frac{t}{N}})\right]= e^{-\frac{nt}{2}}\sum_{k,d=0}^{+\infty} \frac{(-1)^{k}t^k}{k!N^{2d}} \sum_{|\sigma'|=|\sigma|-k+2d}\# \Pi_{k}(\sigma \to \sigma')\; 
p_{\sigma'}(M_{1},\ldots,M_{n}). \label{main eqn M}
\end{equation}
In particular, if $m_{1},\ldots,m_{r}$ denote the lengths of the cycles of $\sigma$, then
\begin{equation}\label{main eqn}
\E[\tr_{N}(B_{\frac{t}{N}}^{m_{1}})\ldots \tr_{N}(B_{\frac{t}{N}}^{m_{r}})] =  e^{-\frac{nt}{2}}
\sum_{k,d=0}^{+\infty}\frac{(-1)^k t^k}{k!N^{2d}}S(\sigma,k,d).
\end{equation}
For all $T\geq 0$, both expansions converge uniformly on $(N,t)\in \N^* \times [0, T]$.
\end{theorem}

In order to understand the role of the defect of a path in our problem, let us write down the result
corresponding to Proposition \ref{lap newt st} for
the functions $p_{\sigma}$. As explained earlier, the number
of cycles of $\sigma \tau$ can be either $\ell(\sigma)+1$ or $\ell(\sigma)-1$, respectively when the
two points exchanged by $\tau$ belong to the same cycle of $\sigma$ or
to two distinct cycles. For each permutation $\sigma\in \S_{n}$, we are led to partition $T_{n}$
into two classes $F(\sigma)$ and $C(\sigma)$, those which fragment a cycle of $\sigma$ and those
which coagulate two cycles. More precisely,
$$F(\sigma)=\{\tau \in T(n): \ell(\sigma\tau)=\ell(\sigma)+1\} \; {\text{and} }\; C(\sigma)= \{\tau
\in T(n): \ell(\sigma\tau)=\ell(\sigma)-1\}.$$
The following result is now a straightforward consequence of Proposition
\ref{lap newt st}.

\begin{proposition} \label{lap newt} Let $\sigma$ be a permutation in $\S_{n}$. Let
$M_{1},\ldots,M_{n}$ be $N\times N$ matrices. Then the following relation holds: 
\begin{eqnarray*}
\frac{1}{2N}\Delta_{\Un} p_{\sigma}(M_{1}U,\ldots,M_{n}U)&=&- \frac{n}{2}\;
p_{\sigma}(M_{1}U,\ldots,M_{n}U) \\
 && \hskip -0.5cm +\sum_{\tau \in F(\sigma)} p_{\sigma\tau}(M_{1}U,\ldots,M_{n}U)+  \frac{1}{N^2}
\sum_{\tau\in C(\sigma)} p_{\sigma\tau}(M_{1}U,\ldots,M_{n}U).
\end{eqnarray*}
\end{proposition}

According to this result, each step which increases the distance to
the identity is penalized by a weight $N^{-2}$. In the proof of the power series expansion, we use
the following lemma.

\begin{lemma} \label{MM} Let $t\geq 0$ and $N>0$ be real numbers. For all $\sigma,\sigma'\in \S_{n}$
and $\epsilon\in \{-1,1\}$, define
$$M^{\epsilon}_{\sigma,\sigma'}=\sum_{k=0}^{+\infty}\frac{\epsilon^k t^{k}}{k!}\frac{\#\Pi_{k}(\sigma\to
\sigma')}{N^{k-(\ell(\sigma')-\ell(\sigma))}}.$$
Then the matrices $(M^{1}_{\sigma,\sigma'})_{\sigma,\sigma'\in \S_{n}}$ and
$(M^{-1}_{\sigma,\sigma'})_{\sigma,\sigma'\in \S_{n}}$ are each other's inverse.
\end{lemma}

\pf Let us define an endomorphism $L$ of $\C[\S_{n}]$ by setting, for all $f\in \C[\S_{n}]$,
$$(Lf)(\sigma)=\sum_{\tau \in F(\sigma)}f(\sigma\tau)+\frac{1}{N^2}\sum_{\tau\in
C(\sigma)}f(\sigma\tau).$$
One checks easily that the matrix $M^{\epsilon}_{\sigma,\sigma'}$ is the matrix of the operator
$e^{\epsilon tL}$ on $\C[\S_{n}]$ and the result follows. \qed


{\noindent \textbf{Proof of Theorem \ref{main expansion} -- }} Consider $T\geq 0$. We claim that the
right-hand side of (\ref{main eqn M}) is a normally convergent series on $(N,t)\in \N^*\times [0,T]$. Indeed, let us define $K=\max\{|p_\sigma(M_1,\ldots,M_n)| : \sigma\in\S_n\}$. Then, for all  $N\geq 1$ and all $t\in[0,T]$, the sum of the absolute values of the terms of the series is smaller than
$$\hskip -.3cm K e^{-\frac{nt}{2}} \sum_{k=0}^{+\infty} \frac{T^k}{k!}\sum_{d=0}^{+\infty}S(\sigma,k,d)= K e^{\frac{n(n-2)}{2}T}.$$
The assertion on the uniform convergence of the expansions follows.

In order to prove (\ref{main eqn M}), we start from the expression given by Proposition \ref{BM
SRW}, at time $\frac{t}{N}$ and with an
arbitrary deterministic initial condition $\pi_{0}=\sigma$. It reads
\begin{equation}\label{pf 1}
\forall \sigma\in \S_{n}, \E\left[\left.
p^{st}_{\pi_{\frac{t}{N}}}(M_{1}B_{\frac{t}{N}},\ldots,M_{n}B_{\frac{t}{N}})\right|\pi_{0}=\sigma\right]= 
e^{-\frac{nt}{2}-\frac{n(n-1)t}{2N}}p^{st}_{\sigma}(M_{1},\ldots,M_{n}).
\end{equation}
We expand the left hand side by using the properties of $(\pi_{t})_{t\geq 0}$. This chain jumps at
rate $\binom{n}{2}$ and its jump chain is a standard discrete-time random walk on the Cayley 
graph of $\S_{n}$, independent of the jump times. Thus, the left-hand side of (\ref{pf 1}) is equal
to
$$\sum_{k=0}^{\infty} e^{-{{n}\choose{2}}\frac{t}{N}}
\binom{n}{2}^{k}\frac{t^k}{k!N^k}\frac{1}{\binom{n}{2}^{k}}\sum_{\sigma'\in
\S_{n}}\sum_{\gamma\in
\Pi_{k}(\sigma\to \sigma')} \E\left[p^{st}_{\sigma'}(M_{1}B_{\frac{t}{N}},\ldots,M_{n}B_{\frac{t}{N}})\right],$$
where the expectation is now only with respect to the Brownian motion. After simplification and switching to normalized traces, (\ref{pf 1}) becomes
$$\hskip -.5cm \forall \sigma\in \S_{n} , \sum_{\sigma'\in
\S_{n}}\E\left[p_{\sigma'}(M_{1}B_{\frac{t}{N}},\ldots,M_{n}B_{\frac{t}{N}})\right]
\sum_{k=0}^{\infty}
\frac{t^k}{k!}\frac{\#\Pi_{k}(\sigma\to\sigma')}{N^{k-(\ell(\sigma')-\ell(\sigma))}}=e^{-\frac{nt}{2}}
p_{\sigma}(M_{1},\ldots,M_{n}).$$
We recognize the expression of $M^1_{\sigma,\sigma'}$ and, by Lemma \ref{MM}, we conclude that for
all $\sigma\in \S_{n}$,
$$ \E\left[p_{\sigma}(M_{1}B_{\frac{t}{N}},\ldots,M_{n}B_{\frac{t}{N}})\right]=e^{-\frac{nt}{2}}
\sum_{k=0}^{\infty}
\frac{(-1)^{k}t^{k}}{k!}\sum_{\sigma'\in \S_{n}}\frac{\#\Pi_{k}(\sigma \to
\sigma')}{N^{k-(\ell(\sigma')-\ell(\sigma))}}p_{\sigma'}(M_{1},\ldots,M_{n}).$$
The first formula follows from the fact that $|\sigma|=n-\ell(\sigma)$. Setting
$M_{1},\ldots,M_{N}$ equal to $I_{N}$ yields the second formula. \qed

Most of the coefficients which appear in the expansion  (\ref{main eqn}) are zero. More precisely,
the situation is the following.

\begin{lemma} \label{dk} Let $\gamma$ be a path of length $k$ and defect $d$ starting at $\sigma$.
Then the following inequalities hold: 
$$0\leq d\leq k \;\;\; {\text{and}} \;\;\; 2d-(\ell(\sigma)-1)\leq k \leq 2d + (n-\ell(\sigma)).$$
In particular, $|k-2d|\leq n-1$. 

Moreover, let $d\geq 0$ be given.
Then $S(\sigma,2d+(n-\ell(\sigma)),d)>0$
and, if $d\geq \ell(\sigma)-1$, then $S(\sigma,2d-(\ell(\sigma)-1),d)>0$. Finally, if
$d\leq \ell(\sigma)-1$, then $S(\sigma,d,d)>0$. 
\end{lemma}

\pf Assume that the path finishes at $\sigma_{k}$. Then the first two inequalities reflect simply the
fact that $1\leq \ell(\sigma_{k})\leq n$.

To prove the second part of the statement, consider $d\geq 0$. Recall that $\sigma$ is fixed. Let
us construct a longest possible path starting at $\sigma$ with defect $d$. For this, we minimize
the defect at each step. First, we build a path by going from $\sigma$ down to the identity through a
geodesic. This takes $n-\ell(\sigma)$ steps and the defect of the path is still zero. Then the path
must make one step up. Immediately after this, it can go down to the identity again. It can repeat
this at most $d$ times without its defect becoming larger than $d$. By then it has length
$2d+(n-\ell(\sigma))$. Thus we have constructed a path of length $2d+(n-\ell(\sigma))$ with defect
$d$. A similar argument works for a shortest path of given defect. \qed 

In particular, for all $d\geq 0$, the contribution of order $N^{-2d}$ to $\E[p_{\sigma}(B_{\frac{t}{N}})]$
 is a polynomial function of $t$ of degree $2d+(n-\ell(\sigma))$ and in which the smallest exponent of $t$ is
$\max(d,2d-(\ell(\sigma)-1))$.

\subsection{Examples} Let us work out explicitly a few examples.

For $n=1$: there is a single path in the Cayley graph of $\S_{1}$. It has length and defect 0.
Thus, we recover the well-known formula
$$\E\left[\tr_{N}(B_{\frac{t}{N}})\right]=e^{-\frac{t}{2}}.$$
For $n=2$: for each $k\geq 0$ there is a unique path of length $k$ starting at the identity. It
has defect $\lfloor \frac{k+1}{2}\rfloor$. Thus,
$$\E\left[\tr_{N}(B_{\frac{t}{N}})^2\right]=e^{-t}\sum_{k=0}^{\infty}\frac{(-1)^{k}t^k}{k! N^{2\lfloor
\frac{k+1}{2}\rfloor}}=e^{-t}\left(\cosh\frac{t}{N}-\frac{1}{N}\sinh\frac{t}{N}\right).$$ 
Similarly, for each $k\geq 0$, there is a unique path of length $k$ starting at $(12)$. It has
defect $\lfloor \frac{k}{2}\rfloor$. Thus,
$$\E\left[\tr_{N}(B_{\frac{t}{N}}^2)\right]=e^{-t}\sum_{k=0}^{\infty}\frac{(-1)^{k}t^k}{k! N^{2\lfloor
\frac{k}{2}\rfloor}}=e^{-t}\left(\cosh \frac{t}{N}-N \sinh\frac{t}{N}\right).$$
For $n=3$: the situation is a bit more complicated but it is still possible to compute everything
by hand. For a path starting at the identity of length $k$ and defect $d$, we must have $2d-2\leq
k\leq 2d$. Hence, if $k$ is odd, it must be equal to $2d-1$. So, for all $l\geq 1$, $S(\id,
2l-1,l)=3^{2l-1}$. If $k$ is even, then two situations are possible. We leave it as an exercise
to check that, for all $l\geq 1$, $S(\id, 2l,l)=3^{2l-1}$ and $S(\id, 2l,l+1)=2.3^{2l-1}$. Finally,
$S(\id,0,0)=1$. We find
$$\E\left[\tr_{N}(B_{\frac{t}{N}})^3\right]=e^{-\frac{3t}{2}}\left(1+\frac{N^2+2}{3N^2}\left(\cosh\frac{3t}{N}-1\right)-
\frac{1}{N}\sinh\frac{3t}{N}\right).$$ 
Similarly, we find
$$\E\left[\tr_{N}(B_{\frac{t}{N}}^2)\tr_{N}(B_{\frac{t}{N}})\right]=e^{-\frac{3t}{2}}\left(\cosh\frac{3t}{N}-
\frac{N^2+2}{3N}\sinh\frac{3t}{N}\right),$$
$$\E\left[\tr_{N}(B_{\frac{t}{N}}^3)\right]=e^{-\frac{3t}{2}}\left(1+\frac{N^2+2}{3}
\left(\cosh\frac{3t}{N}-1\right)-N\sinh\frac{3t}{N}\right).$$
For $n\geq 4$, it seems difficult to determine all the coefficients at once and by hand.
Nevertheless, the following diagram, which indicates how many edges join the various conjugacy
classes of $\S_{4}$ in the Cayley graph allows one to compute specific values of $S(\sigma,k,d)$.
\begin{figure}[h!]
\begin{center}
\scalebox{0.8}{
$$\xymatrix{ & & (123) \ar@<1ex>[dr]^{3} \ar@<1ex>[dl]^{3}  &  \\
     {\id} \ar@<1ex>[r]^{6} & (12) \ar[l]^{1} \ar[ur]^{4} \ar[dr]^{1} &   &  (1234) \ar[ul]^{4} \ar[dl]^{2}
\\
   &  &  (12)(34) \ar@<1ex>[ul]^{2} \ar@<1ex>[ur]^{4} & }$$
}
\end{center}
\caption{The Cayley graph of $\S_{4}$ modulo conjugation.}
\end{figure}
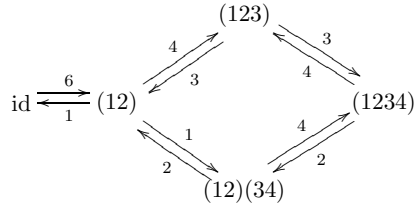

For instance, one can use it to prove the following formulae :
$$e^{2t} \E\left[\tr_{N}(B_{\frac{t}{N}})^{4}\right]= 1+
\frac{1}{N^{2}}(-6t+3t^{2})+\frac{1}{N^{4}}(15t^{2}-20 t^{3} + 5 t^{4}) + O(\frac{1}{N^{6}}) ,$$
$$e^{2t} \E\left[\tr_{N}(B_{\frac{t}{N}}^{4})\right]=(1-6t+8t^2-\frac{8}{3}
t^{3})+\frac{1}{N^{2}}(10 t^{2}-\frac{58}{3}t^{3}+\frac{71}{4}t^{4}-\frac{16}{3}t^{5})+ O(\frac{1}{N^{4}}).$$

\subsection{The case of $SU(N)$}

Let us conclude this section by stating without proof the following analogue of Theorem \ref{main
expansion} in the case of the the special unitary group. This theorem is proved exactly like
its unitary version, by using the observation made at the end of Section \ref{unitary}.

\begin{theorem} \label{main thm su}
Let $N, n\geq 1$ be integers. Let $(B_{t})_{t\geq 0}$ be a
Brownian motion on $SU(N)$ starting at the identity and corresponding to the scalar product $(X,Y)\mapsto -\tr(XY)$ on
$\su(N)$. Let $M_{1},\ldots,M_{n}$ be arbitrary $N\times N$ complex matrices. Let $\sigma$ be an
element of $\S_{n}$. Then, for all $t\geq 0$, we have the following series expansion:
\begin{eqnarray} 
\E\left[p_{\sigma}(M_{1}B_{\frac{t}{N}},\ldots,M_{n}B_{\frac{t}{N}})\right]&=&
e^{-\frac{nt}{2}+\frac{n^2 t}{2N^2}}\sum_{k,d=0}^{+\infty} \frac{(-1)^{k}t^k}{k!N^{2d}}\times \nonumber
\\
&&\hskip -0.6cm\sum_{|\sigma'|=|\sigma|-k+2d}\# \Pi_{k}(\sigma \to \sigma')\;  p_{\sigma'}(M_{1},\ldots,M_{n}).
\end{eqnarray}
In particular, if $m_{1},\ldots,m_{r}$ denote the lengths of the cycles of $\sigma$, then
\begin{equation}\label{main expansion su}
\E[\tr_{N}(B_{\frac{t}{N}}^{m_{1}})\ldots \tr_{N}(B_{\frac{t}{N}}^{m_{r}})] =  e^{-\frac{nt}{2}+\frac{n^2 t}{2N^2}}
\sum_{k,d=0}^{+\infty}\frac{(-1)^k t^k}{k!N^{2d}}S(\sigma,k,d).
\end{equation}
For all $T\geq 0$, both expansions converge uniformly on $(N,t)\in \N^* \times [0, T]$.
\end{theorem}

\section{A representation-theoretic derivation of the power series expansions}
In this section, we give an alternative derivation of the expansions (\ref{main eqn}) and
(\ref{main expansion su}),
based on the representation theory of the unitary and symmetric groups and the relations between
symmetric functions. This approach is less elementary than the one adopted in the previous
sections but we believe that it is more likely to
allow generalizations. In Section \ref{calcul snkd}, we will use it to
compute some of the coefficients $S(\sigma,k,d)$.

\subsection{Expansion for the unitary group}
The integers $N,n\geq 1$ are fixed throughout this section. We write $\lambda \vdash n$ if $\lambda=(\lambda_{1}\geq
\ldots \geq\lambda_{r}>0)$ is a partition of $n$. The integer $r$ is called
the length of $\lambda$ and we denote it by $\ell(\lambda)$. We denote the set of all partitions
by $\PP$.

Let $\lambda$ be a partition of $n$. We denote by $s_{\lambda}$ the Schur function associated to the
partition $\lambda$, whose definition is given in \cite[I.3]{Macdonald}. For all $U\in U(N)$, the
number $s_{\lambda}(U)$ is defined as the value of $s_{\lambda}$ on the eigenvalues of $U$. We
will use the fact that, if $\ell(\lambda)>N$, then the symmetric polynomial in $N$ variables
determined by $s_{\lambda}$ is the zero polynomial. This follows for example from the expression
of $s_{\lambda}$ as a determinant in the elementary symmetric functions \cite[I.3,
(3.5)]{Macdonald}.

Recall the definition of the power sums, that is, the functions $p^{st}_{\sigma}:U(N)\lra \C$ for 
$\sigma \in \S_{n}$ (see Definition \ref{power sums 1}). 

The Schur functions and the power sums are related as follows. Let $\chi^{\lambda}:\S_{n}\lra
\C$ denote the character of the irreducible representation of $\S_{n}$ associated with $\lambda$.
Then one has the following pair of relations \cite[I.7,(7.7)]{Macdonald}:
\begin{equation}\label{p s}
\forall \lambda \vdash n, \; s_{\lambda}=\frac{1}{n!}\sum_{\sigma\in \S_{n}} \chi^{\lambda}(\sigma)
p^{st}_{\sigma},
\end{equation}
\begin{equation}\label{s p}
\forall \sigma \in \S_{n}, \; p^{st}_{\sigma}=\sum_{\lambda \vdash n} \chi^{\lambda}(\sigma)
s_{\lambda}.
\end{equation}

The set $\widehat{U(N)}$ of isomorphism classes of irreducible representations (irreps) of $U(N)$ is
in one-to-one correspondence with the set $\Z^N_{\downarrow}$ of non-increasing sequences $\alpha=(\alpha_{1}\geq
\ldots \geq \alpha_{N})$ of elements of $\Z$. Even when some of the $\alpha_{i}$'s are negative, the
Schur function $s_{\alpha}$ is well-defined and the character of the irrep $\alpha$ is
$\chi_{\alpha}(U)=s_{\alpha}(U)$.

Let $(B_{t})_{t\geq 0}$ be  the Brownian motion on $U(N)$ of Theorem \ref{main expansion}. Let
$dU$ denote the normalized Haar measure on $U(N)$. For each $t>0$, let $Q_{t}$ denote the heat
kernel at time $t$ on $U(N)$, that is, the density of the distribution of $B_{t}$ with respect to
the Haar measure. Our main result is the following reformulation of (\ref{main eqn}).

\begin{theorem} \label{main 2}
Let $N, n\geq 1$ be integers. Let $\sigma$ be an element of $\S_{n}$. Then, for all $t\geq 0$, 
\begin{equation}\label{main eqn 2}
N^{-\ell(\sigma)}\int_{U(N)} p^{st}_{\sigma}(U) Q_{\frac{t}{N}}(U)\; dU =  e^{-\frac{nt}{2}} \sum_{k,d=0}^{+\infty}\frac{(-1)^k t^k}{k!N^{2d}}S(\sigma,k,d).
\end{equation}
\end{theorem}

In the course of the proof, we admit two lemmas which we prove afterwards. We have preferred this
order to the strict logical order to make the proof easier to follow.\\

\pf If $t=0$, the result is clearly true. When $t>0$, the proof consists in expanding $Q_{t}$ into
the sum of its Fourier series and turning all
quantities related to $U(N)$ into quantities related to $\S_{n}$. For all $t>0$, the function
$Q_{t}$ is smooth on $U(N)$ and invariant by conjugation. It admits the following uniformly
convergent Fourier expansion \cite[Thm 4.4]{Liao}:
\begin{equation}\label{fourier Q}
Q_{t}(U)=\sum_{\alpha\in \Z^N_{\downarrow}}
e^{-\frac{c_{2}(\alpha) t}{2}}s_{\alpha}(I_{N})\overline{s_{\alpha}(U)},
\end{equation}
where the number $c_{2}(\alpha)$ is characterized by the equality
$\Delta_{U(N)}\chi_{\alpha}=-c_{2}(\alpha) \chi_{\alpha}$. Using the relation (\ref{s p}) to
expand $p^{st}_{\sigma}(U)$, we find the following expression for the left-hand side of
(\ref{main eqn 2}):
$$\int_{U(N)} p^{st}_{\sigma}(U)Q_{\frac{t}{N}}(U)\; dU = \sum_{\alpha\in
\Z^N_{\downarrow},\mu\vdash n}
e^{-\frac{c_{2}(\alpha)t}{2N}} s_{\alpha}(I_{N})\chi^{\mu}(\sigma) \int_{U(N)} 
\overline{s_{\alpha}(U)}s_{\mu}(U)\;
dU.$$
By the orthogonality properties of the characters of irreps, the integral in the
right-hand side is zero unless $\alpha=\mu$. Hence, we can replace the sum over $\alpha$ and $\mu$
by a sum over the partitions $\mu$ such that $\mu\vdash n$ and $\ell(\mu)\leq N$:
\begin{equation}\label{integral inter 1}
\int_{U(N)} p^{st}_{\sigma}(U)Q_{\frac{t}{N}}(U)\; dU = \sum_{\mu\vdash n,\ell(\mu)\leq N}
e^{-\frac{c_{2}(\mu)t}{2N}} s_{\mu}(I_{N})\chi^{\mu}(\sigma).
\end{equation}
We still need to express $s_{\mu}(I_{N})$ and $c_{2}(\mu)$ in terms of
quantities related to the symmetric group.

In order to compute $s_{\mu}(I_{N})$, we use the relation (\ref{p s}). Let us define $\Omega=\sum_{\sigma \in \S_{n}}
N^{\ell(\sigma)}\sigma$. This notation is borrowed from \cite{GrossTaylor2}. Then (\ref{p s}) implies
the equality
\begin{equation}
\label{D 1}
s_{\mu}(I_{N})=\frac{1}{n!}\chi^{\mu}(\Omega).
\end{equation}
In Lemma \ref{Phi}, we will prove that $\chi^{\mu}(\Omega)=0$ if $\ell(\mu)>N$. This allows
us to drop the restriction $\ell(\mu)\leq N$ in the summation. 

Let us compute $c_{2}(\mu)$, the eigenvalue of $\Delta_{U(N)}$ associated to $s_{\mu}$. Thanks to
(\ref{delta p st}), we know the value of $\Delta_{U(N)}p^{st}_{\sigma}$ for all $\sigma$ and (\ref{p s}) expresses
$s_{\mu}$ as a linear combination of power sums. Combining these two equations, we find
\begin{equation}\label{inter 1}
\Delta_{U(N)}s_{\mu}=-Nns_{\mu}-\frac{2}{n!}\sum_{\sigma \in \S_{n}}
\sum_{\tau \in T_{n}} \chi^{\mu}(\sigma \tau) p_{\sigma}^{st},
\end{equation}
where $T_{n}$ is the set of the transpositions of $\S_{n}$. We now use the
following consequence of Schur's lemma: whenever $x$
belongs to the group algebra $\C[\S_{n}]$ and $y$ to the centre of the group algebra,
\begin{equation}
\label{orth}
\forall \mu \vdash n, \;\; \chi^{\mu}(xy)=\frac{\chi^{\mu}(x)\chi^{\mu}(y)}{\chi^{\mu}(1)}.
\end{equation}
This relation implies that the last term of (\ref{inter 1}) is equal to $n(n-1)
\frac{\chi^{\mu}((12))}{\chi^{\mu}(1)}s_{\mu}$. Hence,
\begin{equation}\label{c2 mu}
c_{2}(\mu)=nN +n(n-1)\frac{\chi^{\mu}((12))}{\chi^{\mu}(1)}.
\end{equation}
Combining (\ref{integral inter 1}), (\ref{D 1}) and (\ref{c2 mu}), we find
$$\int_{U(N)} p^{st}_{\sigma}(U)Q_{\frac{t}{N}}(U)\; dU = e^{-\frac{nt}{2}} \sum_{k\geq
0}\frac{(-t)^k}{k!} \sum_{\mu\vdash n}
\frac{\chi^{\mu}(\Omega)\chi^{\mu}(\sigma)}{n!}\left(\frac{n(n-1)}{2N}\frac{\chi^{\mu}((12))}{\chi^{\mu}(1)}\right)^k.$$
By Lemma \ref{compute S} below, the sum over $\mu$ is equal to $N^{\ell(\sigma)}\sum_{d\geq 0}
N^{-2d}S(\sigma,k,d)$. The result follows immediately. \qed

\begin{lemma} \label{compute S} 
Let $\sigma$ be an element of $\S_{n}$. Let $k\geq 0$ be an integer.
Then
\begin{equation}
\label{S chi}
\sum_{d\geq 0} \frac{S(\sigma,k,d)}{N^{2d}} = N^{-\ell(\sigma)-k}\sum_{\mu\vdash n}
\frac{\chi^{\mu}(\Omega)\chi^{\mu}(\sigma)}{n!}
\left(\frac{n(n-1)}{2}\frac{\chi^{\mu}((12))}{\chi^{\mu}(1)}\right)^k.
\end{equation}
\end{lemma}

\pf There is no issue of convergence, since the sum on the left-hand side is finite. Let
$\gamma=(\sigma_{0},\ldots,\sigma_{k})$ be a path of defect $d$. By Lemma \ref{l d k},
$\ell(\sigma_{k})=\ell(\sigma_{0})+k-2d$. Hence,
$$\sum_{d\geq 0} \frac{S(\sigma,k,d)}{N^{2d}} = N^{-\ell(\sigma)-k}\sum_{\sigma'\in \S_{n}}
\#\Pi_{k}(\sigma \to \sigma') N^{\ell(\sigma')}.$$
Now, $\#\Pi_{k}(\sigma \to \sigma')$ is the number of $k$-tuples $(\tau_{1},\ldots,\tau_{k})\in
T_{n}^k$ such that $\sigma\tau_{1}\ldots \tau_{k}=\sigma'$. A standard computation based on the 
fact that $\sum_{\mu\vdash n}\chi^{\mu}(1)\chi^{\mu}(\sigma)=n! \delta_{\sigma,\id}$
and on (\ref{orth}) leads to 
$$\#\Pi_{k}(\sigma \to \sigma')=\sum_{\mu\vdash n}
\frac{\chi^{\mu}(1)\chi^{\mu}(\sigma^{-1}\sigma')}{n!}
\left(\frac{n(n-1)}{2}\frac{\chi^{\mu}((12))}{\chi^{\mu}(1)}\right)^k.$$
The result follows by summing over $\sigma'$ and applying (\ref{orth}) again. \qed

\begin{lemma} \label{Phi}
Define $\Omega \in \C[\S_{n}]$ by $\Omega =\sum_{\sigma \in \S_{n}}
N^{\ell(\sigma)}\sigma$. Then the following relations hold.\\
1.  For all $\mu\vdash n$ such that $\ell(\mu)\leq N$, $\chi^{\mu}(\Omega)=n! s_{\mu}(I_{N})$.\\
2. For all $\mu\vdash n$ such that $\ell(\mu)> N$, $\chi^{\mu}(\Omega)=0$.
\end{lemma}

\pf The first assertion follows immediately from (\ref{p s}).

In order to prove the second assertion, let us introduce the Jucys-Murphy elements
$X_{1},\ldots,X_{n}$ of $\C[\S_{n}]$, defined by $X_{1}=0$ and $X_{i}=(1\; i)+(2\; i)+\ldots + (i-1\;i)$
for $i\in \{2,\ldots, n\}$. They generate a maximal Abelian subalgebra of $\C[\S_{n}]$. In
particular, they can be simultaneously diagonalized in every irreducible representation of
$\S_{n}$. We borrow the following statements from \cite{VershikOkounkov}.

Let $\mu=(\mu_{1},\ldots,\mu_{r})$ be a partition of $n$. The subset
$D_{\mu}=\{(i,j)\in (\N^*)^2 : j\leq \mu_{i}\}$ of $\Z^2$ is called the diagram of
$\mu$. An element $(i,j)$ of $D_{\mu}$ is called a box and its content is defined as the
integer $c(i,j)=j-i$. 

The space of the irreducible representation of $\S_{n}$ associated to $\mu$ admits a basis
which diagonalizes the Jucys-Murphy elements and is indexed by the bijections $t:D_{\mu}\lra
\{1,\ldots,n\}$ which are increasing in each variable. These bijections are usually called tableaux.
The eigenvalue of the Jucys-Murphy element $X_{k}$ on the vector associated to the tableau $t$ is the
content of the box $t^{-1}(k)$. 

We need also the following well-known fact: for every $k \in \{0,\ldots,n-1\}$, the $k$-th elementary
symmetric function of the Jucys-Murphy elements is equal to $\sum_{\sigma \in \S_{n}} {\bm 1}_{|\sigma|=k}
\; \sigma$, where $|\sigma|=n-\ell(\sigma)$. This can be proved as follows. For each $m\in
\{1,\ldots,n\}$, let us imbed $\S_{m}$ into $\S_{n}$ as the subgroup which leaves $\{m+1,\ldots,n\}$
invariant. For all $m\in \{1,\ldots,n\}$ and $k\in\{0,\ldots,m-1\}$, set $\Sigma_{m,k}=\sum_{\sigma\in
\S_{m}} {\bm 1}_{|\sigma|=k} \; \sigma$. Let us also define $e_{m,k}=\sum_{i_{1}<\ldots <i_{k}\leq
m} X_{i_{1}}\ldots X_{i_{k}}$. We need to prove that $e_{m,k}=\Sigma_{m,k}$. This is clearly true if
$k\in\{0,1\}$. The general case follows by induction on $m$, each inductive step being proved by
induction on $k$, thanks to the relations
$$\Sigma_{m,k}=\Sigma_{m-1,k}+\Sigma_{m-1,k-1} X_{k} \;\;\; \mbox{ and }\;\;\; e_{m,k}=e_{m-1,k}+e_{m-1,k-1}
X_{k}.$$

From the equality proved in the last paragraph and the relation $|\sigma|=n-\ell(\sigma)$, we
deduce the following equality in the polynomial ring $\C[\S_{n}][z]$: 
\begin{equation}\label{id Sn poly}
\prod_{i=1}^{n} (z+X_{i})=\sum_{\sigma\in \S_{n}} z^{\ell(\sigma)} \sigma.
\end{equation}
Evaluating at $z=N$ and applying $\chi^{\mu}$, we find
$$\chi^{\mu}(\Omega)=\chi^{\mu}\left(\prod_{i=1}^{n} (N+X_{i})\right)=\sum_{t \mbox{ \rm
\scriptsize tableau}}\prod_{i=1}^{n}
(N+c(t^{-1}(i)))=\chi^{\mu}(1) \prod_{i=1}^{\ell(\mu)}\prod_{j=1}^{\mu_{i}} (N+j-i).$$

If $\ell(\mu)\geq N+1$, then $(N+1,1)$ is a box of $D_{\mu}$, whose content is $-N$. It follows
that $\chi^{\mu}(\Omega)=0$ in this case. \qed

\subsection{Expansion for the special unitary group}

Let us apply a similar analysis to the special unitary group in order to derive (\ref{main
expansion su}) in another way. 

By restriction, any irrep of $U(N)$ determines an irrep of $SU(N)$ and the restrictions of
$\alpha,\alpha'\in \Z^N_{\downarrow}$ are isomorphic if and only if there exists $k\in \Z$ such that
$\alpha'=(\alpha_{1}+k,\ldots,\alpha_{N}+k)$. Hence, the set of irreps of $SU(N)$ is in
one-to-one correspondence with the set of partitions of length at most $N-1$ and the character of
the irreducible representation corresponding to a partition $\lambda$ is given by the Schur function
$s_{\lambda}$.

Let $(B_{t})_{t\geq 0}$ be  the Brownian motion on $SU(N)$ of Theorem \ref{main expansion su}. Let
$dU$ denote the Haar measure on $SU(N)$. For each $t>0$, let $Q_{t}$ denote the heat kernel at time
$t$ on $SU(N)$, that is, the density of the distribution of $B_{t}$ with respect to the Haar
measure. We reformulate (\ref{main expansion su}) as follows.

\begin{theorem} \label{main 2 su}
Let $N, n\geq 1$ be integers. Let $\sigma$ be an element of $\S_{n}$. Then, for all $t\geq 0$, 
\begin{equation}
N^{-\ell(\sigma)}\int_{SU(N)} p^{st}_{\sigma}(U) Q_{\frac{t}{N}}(U)\; dU = 
e^{-\frac{nt}{2}+\frac{tn^2}{2N^2}} \sum_{k,d=0}^{+\infty}\frac{(-1)^k t^k}{k!N^{2d}}S(\sigma,k,d).
\end{equation}
\end{theorem}

\pf If $t=0$, both sides are equal to 1. Assume that $t>0$. The Fourier expansion of $Q_{t}$ is
then the following:
\begin{equation}
\label{fourier Q su}
Q_{t}(U)=\sum_{\substack{\lambda\in \PP \\  \ell(\lambda)\leq N-1}}
e^{-\frac{c'_{2}(\lambda) t}{2}}s_{\lambda}(I_{N})\overline{s_{\lambda}(U)},
\end{equation}
where now $c_{2}'(\lambda)$ is defined by the equality
$\Delta_{SU(N)}\chi_{\lambda}=-c'_{2}(\lambda) \chi_{\lambda}$. Combined with the relation (\ref{s
p}), it implies 
$$ \int_{SU(N)} p^{st}_{\sigma}(U)Q_{\frac{t}{N}}(U)\; dU = \sum_{\substack{\lambda\in
\PP,\mu\vdash n \\ \ell(\lambda)\leq N-1}} e^{-\frac{c'_{2}(\lambda)t}{2N}}
s_{\lambda}(I_{N})\chi^{\mu}(\sigma) \int_{SU(N)} 
\overline{s_{\lambda}(U)}s_{\mu}(U)\; dU.$$
When $\ell(\mu)>N$, $s_{\mu}$ is identically zero on $SU(N)$. When $\ell(\mu)\leq N-1$, then the
integral in the right-hand side of the last equation is equal to $\delta_{\lambda,\mu}$. Let us
consider the terms of the sum for which $\ell(\mu)=N$. In this case, let us write
$\mu=\mu'+(\mu_{N},\ldots,\mu_{N})$, so that $\ell(\mu')\leq N-1$. Then the integral is equal to
$\delta_{\lambda,\mu'}$ and we may assume that $\lambda=\mu'$. In this case, $s_{\lambda}(I_{N})=s_{\mu'}(I_{N})=s_{\mu}(I_{N})$
and $c'_{2}(\lambda)=c'_{2}(\mu')$. Hence,
$$
\int_{SU(N)} p^{st}_{\sigma}(U)Q_{\frac{t}{N}}(U)\; dU = \hskip -0.3cm\sum_{\mu\vdash n,\ell(\mu)\leq N-1}
e^{-\frac{c'_{2}(\mu)t}{2N}} s_{\mu}(I_{N})\chi^{\mu}(\sigma)  + \sum_{\mu\vdash n,\ell(\mu)= N}
e^{-\frac{c'_{2}(\mu')t}{2N}} s_{\mu}(I_{N})\chi^{\mu}(\sigma).$$
Let us first compute $c'_{2}(\mu)$ when $\ell(\mu)\leq N-1$. For this, we use the fact that $s_{\mu}$
is an eigenvector of $\Delta_{U(N)}$ whose restriction to $SU(N)$ is $s_{\mu}$. Since
$\Delta_{SU(N)}p_{\sigma}^{st}=\left(\Delta_{U(N)}+\frac{n^2}{N}\right)p_{\sigma}^{st}$ whenever
$\sigma\in \S_{n}$, we find, thanks to (\ref{orth}),
$$\Delta_{SU(N)}s_{\mu}=\left(-Nn-n(n-1)
\frac{\chi^{\mu}((12))}{\chi^{\mu}(1)}+\frac{n^2}{N}\right) s_{\mu}.$$
When $\ell(\mu)=N$, we are interested in $s_{\mu'}$ but $s_{\mu}$ is still an
eigenvector of $\Delta_{U(N)}$ whose restriction to $SU(N)$ is $s_{\mu'}$. Hence,
$$\Delta_{SU(N)}s_{\mu'}=\Delta_{SU(N)}s_{\mu}=\left(-Nn-n(n-1)
\frac{\chi^{\mu}((12))}{\chi^{\mu}(1)}+\frac{n^2}{N}\right) s_{\mu}=-c_{2}'(\mu')s_{\mu'}.$$
Thus, in both sums, the argument of the exponential is
$-\frac{nt}{2}+\frac{n^2t}{2N^2}-t\frac{n(n-1)}{2N}\frac{\chi^{\mu}((12))}{\chi^{\mu}(1)}$. Using
this fact and the first assertion of Lemma \ref{Phi}, we find
$$\int_{SU(N)} p^{st}_{\sigma}(U)Q_{\frac{t}{N}}(U)\; dU = e^{-\frac{nt}{2}+\frac{n^2t}{2N^2}} \sum_{\mu\vdash
n,\ell(\mu)\leq N} e^{-t\frac{n(n-1)}{2N}\frac{\chi^{\mu}((12))}{\chi^{\mu}(1)}}
\frac{\chi^{\mu}(\Omega)\chi^{\mu}(\sigma)}{n!}.$$
The second assertion of Lemma \ref{Phi} tells us that we can remove the restriction $\ell(\mu)\leq
N$, since the other terms are zero. After expanding the exponential, Lemma \ref{compute S} allows us
to finish the proof just as in the unitary case. \qed

\section{Computation of $S((1\ldots n),k,d)$.}
\label{calcul snkd}
In this section, we apply the methods of representation theory to the computation of some of the
coefficients which appear in our main expansions, namely the coefficients $S((1\ldots n),k,d)$ for all $n,k,d\geq 0$.

Let us recall the definition of the Stirling cycle numbers, or Stirling numbers of the first kind
$s(n,k)$, also denoted by $n\brack k$. They are characterized by the identities in $\C[x]$
$$x(x-1)\ldots (x-n+1)=\sum_{k=0}^{n} {n\brack k} x^{k},$$
valid for all $n\geq 0$. In other words,
$${n\brack k}= (-1)^{n-k} e_{n-k}(1,\ldots,n-1)= (-1)^{n-k} \sum_{1\leq i_{1}<\ldots < i_{n-k}\leq
n-1} i_{1}\ldots i_{n-k},$$
where $e_{n-k}$ denotes the $(n-k)$-th elementary symmetric function.


By applying the alternating character to the identity (\ref{id Sn poly}), we find the relation
$$\sum_{\sigma\in \S_{n}} \epsilon(\sigma) x^{\ell(\sigma)}=x(x-1)\ldots (x-n+1)=\sum_{k=0}^{n} {n\brack k}
x^{k},$$
from which we deduce that $\left| {n\brack k}\right|$ is the number of elements of $\S_{n}$ with
exactly $k$ cycles, or in other words at distance $n-k$ from the identity. In particular, ${n\brack
0}=0$. Let us make the convention that ${n\brack k}=0$ if $k<0$. The main result of this section is
the following.

\begin{proposition} \label{snkd}
For all $n,k,d\geq 0$, 
$$S((1\ldots n),k,d)=\frac{1}{n} \sum_{\substack{r,s,l,m\geq 0 \\ r+s=n-1\\ l+m=n-1-k+2d}}
\!\!\frac{(-1)^{l+r}}{r!s!}\left(\frac{n}{2}(s-r)\right)^k {s+1 \brack s+1-l}{r+1
\brack r+1-m}.$$

\end{proposition} 

\pf Instead of computing $S((1\ldots n),k,d)$ we compute the sum of $S(\sigma,k,d)$ when $\sigma$
spans the set of all $n$-cycles. Dividing the result by $(n-1)!$ yields $S((1\ldots n),k,d)$. Now,
a path of length $k$ starting at an $n$-cycle has defect $d$ if and only if it ends at a distance
$n-1-k+2d$ from the identity. Let us recall some of the notation used in the proof of Lemma \ref{Phi}. The integer $n$ being
fixed, we set $\Sigma_{r}=\sum_{|\sigma|=r} \sigma$.  Hence,
\begin{eqnarray}
S((1\ldots n),k,d) &=& \frac{1}{(n-1)!} \sum_{\sigma \in \S_{n},  |\sigma|=n-1}
S(\sigma,k,d)\nonumber \\
&=& \frac{1}{(n-1)!} \sum_{\substack{\sigma,\pi \in \S_{n}, 
|\sigma|=n-1, |\pi|=n-1-k+2d}} \#\Pi_{k}(\sigma \to \pi) \nonumber \\
&=& \frac{1}{(n-1)!} \sum_{\lambda \vdash n} 
\frac{\chi^\lambda(\Sigma_{n-1})\chi^{\lambda}(\Sigma_{n-1-k+2d})}{n!}
\left(\frac{\chi^{\lambda}(\Sigma_{1})}{\chi^{\lambda}(\id)}\right)^{k}. \label{somme}
\end{eqnarray}

Now we use the following fact, which is a consequence of the description of the representations of
$\S_{n}$ given by Okounkov and Vershik \cite{VershikOkounkov} and recalled briefly in the proof of
Lemma \ref{Phi}: 
\begin{equation}\label{car sphere}
\forall r\in\{0,\ldots,n-1\} \; , \;\; \frac{\chi^{\lambda}(\Sigma_{r})}{\chi^{\lambda}(\id)} =
e_{r}\left(\{ c(\square) : \square \in \lambda\}\right).
\end{equation}
In words, the right-hand side of this equation is the $r$-th elementary symmetric function of the
contents of the boxes of the diagram of $\lambda$. In particular, if the diagram of $\lambda$ has
at least two boxes of content $0$, then
$\chi^{\lambda}(\Sigma_{n-1})=0$. Hence, the non-zero terms of the sum (\ref{somme}) arise from
the partitions which are hooks, that is, of the form $\eta_{r}=(n-r \; 1^{r})$ for some $r\in
\{0,\ldots, n-1\}$. This fact is well known (see for example the appendix of \cite{LandoZvonkin}),
and $\eta_{r}$ is the representation $\bigwedge^{r}\St$, where
$\St$ is the restriction of the natural representation of $\S_{n}$ on $\C^n$ to the hyperplane of
equation $\{z_{1}+\ldots+z_{n}=0\}$. This representation is of degree $\binom{n-1}{r}$ and
$\chi^{\eta_{r}}((1\ldots n))=(-1)^{r}$. 

Let us introduce the notation $s=n-1-r$. It follows easily from (\ref{car sphere}) that
\begin{equation}\label{S1}
\forall r\in\{0,\ldots,n-1\} \; , \;\;
\frac{\chi^{\eta_{r}}(\Sigma_{1})}{\chi^{\eta_{r}}(\id)}=\frac{n(n-1)}{2}-nr 
=\frac{n}{2}(s-r).
\end{equation}
 
 The contents of the boxes of $\eta_{r}$ are $\{-r, \ldots, 0,\ldots,s\}$. Hence, by the
definition of the Stirling numbers and (\ref{car sphere}), we have for all $r\in \{0,\ldots,n-1\}$
\begin{eqnarray}
\frac{\chi^{\eta_{r}}(\Sigma_{n-1-k+2d})}{\chi^{\eta_{r}}(\id)} &=&
e_{n-1-k+2d}(-r,\ldots,0,\ldots,s)\nonumber\\
&=& \sum_{\substack{l,m\geq 0\\ l+m=n-1-k+2d}} e_{l}(1,\ldots,s)
(-1)^{m}e_{m}(1,\ldots,r)\nonumber \\
&=& \sum_{\substack{l,m\geq 0\\ l+m=n-1-k+2d}} (-1)^l {s+1 \brack s+1-l}{r+1 \brack
r+1-m}.\label{Snk}
\end{eqnarray}

Finally, combining (\ref{S1}) and (\ref{Snk}), we find
the expected result. \qed

It seems that Proposition \ref{snkd} should allow one to find a simple generating function for the
numbers $S((1\ldots n),k,d)$. Our best result in this direction is the following. We use the
notation $(x)_{n}=x(x-1)\ldots (x-n+1)$. 

\begin{proposition} For all $n,N\geq 0$, one has
\begin{eqnarray*}
\sum_{k,d\geq 0} \frac{(-1)^{k}t^{k}}{d! N^{2d}} S((1\ldots
n),k,d)&=&\frac{(-1)^{n}}{n}\sum_{\substack{r,s\geq 0\\
r+s=n-1}}\frac{1}{r!s!}e^{\frac{1}{4N^2}(s-r)^2 n^2t^2}  \left(\frac{nt}{2} (s-r)\right)_{s+1}\left(-\frac{nt}{2} (s-r)\right)_{r+1}.
\end{eqnarray*}
\end{proposition}

We emphasize that this generating function is, unfortunately, exponential with respect to $d$
instead of $k$.

In Proposition \ref{snkd}, when $k$ takes the largest possible value given $n$ and $d$, namely
$n-1+2d$, then $l$ and $m$ must be equal to $0$ and the identity ${n\brack n}=1$ simplifies
greatly the expression. This leads us to the following corollary.

\begin{corollary}\label{decomp cycles arbitrary}
Let $n\geq 1$ be an integer. For each $p\geq 0$, let $c_{n,p}$ denote the number of distinct ways in
which the cycle $(1\ldots n)\in \S_{n}$ can be written as a product of $p$ transpositions. The
number $c_{n,p}$ is non-zero if and only if $p=n-1+2d$ for some $d\geq 0$. In this case,
$$c_{n,p}=S((1\ldots n),n-1+2d,d)=\frac{n^{p}}{n!}\sum_{r=0}^{n-1} (-1)^{r} \binom{n-1}{r}
\left(\frac{n-1}{2}-r\right)^{p}.$$
For each $n\geq 1$, one has the equality $$\sum_{p\geq 0} c_{n,p} \frac{x^p}{p!}=
\frac{1}{n!} e^{\frac{n(n-1)}{2}x}\left(1-e^{-nx}\right)^{n-1}.$$
In particular, $c_{n,n-1}=n^{n-2}$, $c_{n,n+1}=\frac{1}{24}(n^2-1)n^{n+1}$ and
$$c_{n,n+3}=\frac{1}{5760} (5n-7)(n+3)(n+2)(n^2-1)n^{n+3}.$$
\end{corollary}

\begin{remark} The value of $c_{n,n-1}$ is classical. The sequence $(c_{n,n+1})_{n\geq 1}$ is known
as A060603 in the Online Encyclopedia of Integer Sequences \cite{OLEIS}.
\end{remark}

\section{Asymptotic distribution} 

One of the consequences of Theorem \ref{main expansion} is that the limit as $N$ tends to
infinity of $\E[p_{\sigma}(B_{\frac{t}{N}})]$ exists. Using Lemma \ref{dk}, we get the following
result.

\begin{proposition}\label{asympt dist} Consider $\sigma\in \S_{n}$. The limit of $\E[p_{\sigma}(B_{\frac{t}{N}})]$
as $N$ tends to infinity exists and it is equal to 
$$\lim_{N\to\infty}
\E\left[p_{\sigma}(B_{\frac{t}{N}})\right]=e^{-\frac{nt}{2}}\sum_{k=0}^{|\sigma|}
(-1)^{k}\frac{S(\sigma,k,0)}{k!} t^k.$$
\end{proposition}

Unfortunately, Proposition \ref{snkd} does not seem to lead easily to a simple expression for
$S(\sigma,k,0)$ nor even $S((1\ldots n),k,0)$. In this section, we determine a simple expression of
$S(\sigma,k,0)$ for all $\sigma$ and $k\geq 0$.
For this, we prove a factorization property and use the relation between the metric geometry of
the Cayley graph of $\S_{n}$ and the lattice of non-crossing partitions of the cycle $(1\ldots n)$.
The fact that the two expressions of $S((1\ldots n),k,0)$ given by Propositions \ref{snkd} and
\ref{calc S} agree is not obvious, at least for the author.

\subsection{The factorization property} The factorization property is the following result. It
reduces the problem of the determination of $S(\sigma,k,0)$ to the case where $\sigma$ is a cycle.

\begin{proposition}\label{factorization} Let $m_{1},\ldots,m_{r}$ be positive integers. Then
$$\lim_{N\to \infty}\E\left[\tr_{N}(B_{\frac{t}{N}}^{m_{1}})\ldots
\tr_{N}(B_{\frac{t}{N}}^{m_{r}})\right] = \lim_{N\to
\infty}\E\left[\tr_{N}(B_{\frac{t}{N}}^{m_{1}})\right] \ldots \lim_{N\to
\infty}\E\left[\tr_{N}(B_{\frac{t}{N}}^{m_{r}})\right].$$
More precisely,
\begin{equation}\label{fact O}
\E\left[\tr_{N}(B_{\frac{t}{N}}^{m_{1}})\ldots
\tr_{N}(B_{\frac{t}{N}}^{m_{r}})\right]-\E\left[\tr_{N}(B_{\frac{t}{N}}^{m_{1}})\right] \ldots  
\E\left[\tr_{N}(B_{\frac{t}{N}}^{m_{r}})\right] = O(N^{-2}),
\end{equation}
uniformly in $t$ on bounded intervals.
\end{proposition}

We start by proving the following property of the numbers $S(\sigma,k,0)$. It is in fact
equivalent to the proposition.

\begin{proposition} \label{fact S}Consider $\sigma\in \S_{n}$. Assume that $\sigma=c_{1}\ldots
c_{\ell(\sigma)}$
is the decomposition of $\sigma$ as a product of cycles with disjoint support. Then
\begin{equation}\label{eqn fact S}
\forall k\geq 0 \; , \;\; S(\sigma,k,0)=\sum_{l_{1}+\ldots+l_{\ell(\sigma)}=k}
\frac{k!}{l_{1}!\ldots l_{\ell(\sigma)}!} S(c_{1},l_{1},0)\ldots
S(c_{\ell(\sigma)},l_{\ell(\sigma)},0).
\end{equation}
\end{proposition}

\pf The number $S(\sigma,k,0)$ is the number of paths of length $k$ starting at $\sigma$ and which
at each step move towards a permutation with one more cycle than their current position. As
we already observed several times,© each step of such a path corresponds to the multiplication by
a transposition which exchanges two points which belong to the same cycle of $\sigma$. There is
thus a natural partition of the set of all steps of such a path, according to the cycle of $\sigma$
in which their support is contained. Let us introduce some notation. Let
$(\sigma_{0}=\sigma,\sigma_{1},\ldots,\sigma_{k})$ be a path with defect zero. For each $i\in
\{1,\ldots,k\}$, set $\tau_{i}=\sigma_{i-1}^{-1}\sigma_{i}$. Let $(C_{1},\ldots,C_{\ell(\sigma)})$
be the partition of $\{1,\ldots,k\}$ determined by the fact that $i\in C_{j}$ if and only if the
support of $\tau_{i}$ is contained in the support of $c_{j}$. Then it is clear that, for all $j\in
\{1,\ldots,\ell(\sigma)\}$, the transpositions $(\tau_{i}, i\in C_{j})$ are the steps of a path with
defect zero starting from $c_{j}$.

Hence, constructing a path of length $k$ starting at $\sigma$ and with defect zero is equivalent to
constructing a collection of $\ell(\sigma)$ paths with defect zero starting at $c_{1},\ldots,c_{\ell(\sigma)}$
respectively, whose lengths $l_{1},\ldots,l_{\ell(\sigma)}$ add up to $k$, and a shuffling of the
steps of these paths, that is, a sequence $(C_{1},\ldots,C_{\ell(\sigma)})$ of subsets of $\{1,\ldots,k\}$ which
partition $\{1,\ldots,k\}$ and whose cardinals are $l_{1},\ldots,l_{\ell(\sigma)}$ respectively. 
The equation (\ref{eqn fact S}) is just the translation in symbols of the last sentence. \qed

{\noindent \textbf{Proof of Proposition \ref{factorization} -- }} By Theorem \ref{main expansion} and
Proposition \ref{fact S}, the terms of degree $N^0$ of the difference on the left hand side of (\ref{fact
O}) vanish. Hence, this difference is of the form $N^{-2} F(t,N^{-2})$ for some entire function
$F$. The result follows. \qed 

We have observed after Definition \ref{def S} that $S(\sigma,k,d)$ depends only on the conjugacy
class of $\sigma$. Hence, we need to compute $S((1\ldots m),k,0)$. The arguments
of the proof of Proposition \ref{fact S} show that the paths of defect 0 starting at $(1\ldots m)$
stay in $\S_{m}$ if we identify $\S_{m}$ with the subgroup of $\S_{n}$ which leaves $\{m+1,\ldots,n\}$
invariant. Hence, we are reduced to the computation of $S((1\ldots n),k,0)$ for all $n\geq 1$ and
$k\in \{0,\ldots,n-1\}$. This computations involves non-crossing partitions. For the sake of being
self-contained, we give a brief review of the properties of non-crossing partitions that we
use.

\subsection{Non-crossing partitions} 
Let $P=\{P_{1},\ldots,P_{\ell}\}$ be a partition of $\{1,\ldots,n\}$. The partition $P$ is said to
be non-crossing if there does not exist $i,j,k,l\in\{1,\ldots,n\}$ such that $i<j<k<l$
and $r,s\in\{1,\ldots,\ell\}$ with $r\neq s$ such that $i,k \in P_{r}$ and $j,l \in P_{s}$. 
Another way to formulate the fact that $P$ is non-crossing is the following. For each class $P_{j}$
of the partition, let $H_{j}$ denote the convex hull in $\C$ of $\{e^{\frac{2ik\pi }{n}}: k\in P_{j}\}$.
Then $P$ is non-crossing if and only if for all $i,j\in \{1,\ldots, \ell\}$, $i\neq j \Rightarrow
H_{i}\cap H_{j}=\varnothing$. This notion is relative to the cyclic order on $\{1,\ldots,n\}$
determined by $(1\ldots n)$.  We denote
by $NC(n)$ the set of non-crossing partitions
of the cycle $(1\ldots n)$. This set has been first considered by Kreweras in \cite{Kreweras}.

The fineness relation between partitions restricted to $NC(n)$ makes $NC(n)$ a poset. More
precisely, we say that $P_{1}\preccurlyeq P_{2}$ if every class of $P_{1}$ is contained in a class
of $P_{2}$. The poset $(NC(n),\preccurlyeq)$ is
in fact a lattice, which means that suprema and infima exist. There is in particular a maximum,
$\{\{1,\ldots,n\}\}$, which we denote by $1_{n}$, and a minimum, $\{\{1\},\ldots,\{n\}\}$, which
we denote by $0_{n}$. The
poset $NC(n)$ can be made into a graph by joining two partitions $P$ and $Q$ if they are distinct
and comparable, say $P\prec Q$, and the interval $[P,Q]=\{R\in NC(n): P\preccurlyeq R \preccurlyeq
Q\}$ is reduced to $\{P,Q\}$.

A non-crossing partition $P=(P_{1},\ldots,P_{\ell})$ of the cycle  $(1\ldots n)$ determines an
element $\sigma_{P}$ of $\S_{n}$ as follows: take the cycles of $\sigma_{P}$ to be the classes of $P$ with the
cyclic order induced by $(1\ldots n)$. In symbols, if $i\in P_{j}$, then
$$\sigma_{P}(i)=(1\ldots n)^{k}(i) \; , \;\; \text{where} \;\; k=\min\{l\geq 1: (1\ldots
n)^{l}(i)\in P_{j}\}.$$
In particular, $\sigma_{0_{n}}=\id$ and $\sigma_{1_{n}}=(1\ldots n)$. The partial order on $NC(n)$
corresponds via the mapping $P\mapsto \sigma_{P}$ to the following partial order on $\S_{n}$.

Consider $\sigma_{1},\sigma_{2}\in \S_{n}$. Recall that $|\sigma_{1}|=n-\ell(\sigma_{1})$, the
minimal number of terms of a decomposition of $\sigma_{1}$ in a product of transpositions, is the
distance from $\id$ to $\sigma_{1}$ in the Cayley graph of $\S_{n}$ generated by $T_{n}$. By
definition,  we say that $\sigma_{1}\preccurlyeq\sigma_{2}$ if
$|\sigma_{2}|=|\sigma_{1}|+|\sigma_{1}^{-1}\sigma_{2}|$. In words,
$\sigma_{1}\preccurlyeq\sigma_{2}$ if and only if there exists a geodesic path from $\id$ to
$\sigma_{2}$ through $\sigma_{1}$. The identity is the minimum of $\S_{n}$ for this partial order,
and the $n$-cycles the (pairwise incomparable) maximal elements. The next lemma is well known and 
its proof is left to the reader.

\begin{lemma} The mapping from $NC(n)$ to $\S_{n}$ which sends a partition $P$ to the permutation
$\sigma_{P}$ is an isomorphism of posets from $NC(n)$ onto
$[\id,(1\ldots n)]=\{\sigma \in \S_{n} :  \sigma \preccurlyeq
(1\ldots n)\}.$
\end{lemma}

As a consequence of this Lemma, $S((1\ldots n),k,0)$ is the number of decreasing paths of length $k$
starting at $1_{n}$ in $NC(n)$. It turns out to be easier to count increasing paths in
$NC(n)$ starting at $0_{n}$. They are in one-to-one correspondence by the following duality
property of $NC(n)$ discovered by Kreweras.

For $\sigma\in [\id,(1\ldots n)]$, let us introduce $K(\sigma)=\sigma^{-1}(1\ldots n)$. It is
readily checked that $K$ is a decreasing bijection of $[\id,(1\ldots n)]$. The corresponding
decreasing bijection of $NC(n)$ is called the Kreweras complementation. It can be described
combinatorially at the level of non-crossing partitions as follows. 

Given a partition $R$ of $\{1,\ldots,n\}$ and a partition $S$ of $\{1,\ldots,n\}\simeq \{\bar{1},\ldots,\bar{n}\}$,
let $R\cup S$ denote the partition of $\{1,\bar{1},2,\bar{2},\ldots,n,\bar{n}\}$ obtained by merging
$R$ and $S$. Even if $R$ and $S$ are non-crossing, $R\cup S$ may be crossing with respect to the
cyclic order $(1,\bar{1},2,\bar{2},\ldots,n,\bar{n})$. Now let $P$ be a non-crossing
partition of $\{1,\ldots,n\}$. The partition $K(P)$ is by definition the largest element of
$NC(n)$ such that $P\cup K(P)$ is non-crossing.

\begin{figure}[h!]
\begin{center}
\scalebox{0.65}{\includegraphics{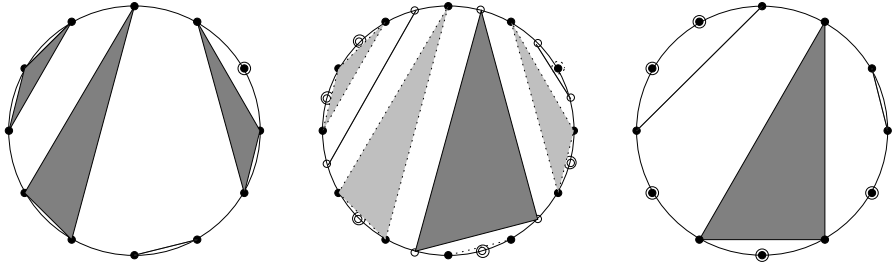}}
\end{center}
\caption{The Kreweras complement of
$\{\{1,3,12\},\{2\},\{4,8,9\},\{5,6,7\},\{10,11\}\}$ is $\{\{1,2\},\{3,9,11\},\{4,7\},\{5\},\{6\},\{8\},\{10\},\{12\}\}$.}
\end{figure}

The following result summarizes this discussion of non-crossing partitions in relation to our problem.

\begin{proposition} For all $n\geq 1$ and $k\geq 0$, $S((1\ldots n),k,0)$ is the number of
increasing paths of length $k$ starting at $\{\{1\},\ldots,\{n\}\}$ in the lattice of non-crossing
partitions of $(1\ldots n)$. 
\end{proposition}

It remains to count these paths. To do this, we use the fact that $S((1\ldots n),n-1,0)=n^{n-2}$.
This is a classical result of combinatorics, since $S((1\ldots n),n-1,0)$ is the number of ways to
write a $n$-cycle as a product of $n-1$ transpositions. It is also a special case of Corollary
\ref{decomp cycles arbitrary}.

\begin{proposition} \label{calc S} For all $n\geq 1$ and $k\geq 0$,
$$S((1\ldots n),k,0)=\binom{n}{k+1} n^{k-1}.$$
It is understood that this number is zero if $k\geq n$.
\end{proposition}

\pf We count the increasing paths of length $k$ in $NC(n)$ starting at $\{\{1\},\ldots,\{n\}\}$ by
first regrouping them according to their terminal point. The possible terminal points of these
paths are exactly the non-crossing partitions of $(1\ldots n)$ into $n-k$ classes. Such partitions
may be classified according to the number of singletons they contain, the number of pairs, and so
on.

Let $s_{1},\ldots,s_{n}$ be non-negative integers such that $s_{1}+\ldots + s_{n}=n-k$ and
$s_{1}+2s_{2}+\ldots + ns_{n}=n$. We say that a partition is of type $(s_{1},\ldots,s_{n})$ if it
contains exactly $s_{i}$ classes of cardinal $i$ for each $i\in \{1,\ldots,n\}$. The number of
non-crossing partitions of type $(s_{1},\ldots,s_{n})$ has been computed by Kreweras
\cite{Kreweras}. It is equal to $\displaystyle \frac{n!}{(k+1)!s_{1}!\ldots s_{n}!}$.

Let $P$ be a non-crossing partition of type $(s_{1},\ldots,s_{n})$. An
increasing path from $\{\{1\},\ldots,\{n\}\}$ to $P$ looked at in the reverse direction is a 
decreasing path from $P$ to $\{\{1\},\ldots,\{n\}\}$. There are as many such paths as there are
geodesic paths from the permutation $\sigma_{P}$ induced by $P$ to the identity, that is,
$S(\sigma_{P},k,0)$. Let us apply Lemma \ref{fact S} to compute this number. The only non-zero term
in the sum corresponds to the situation where $l_{i}=|c_{i}|$ for each $i\in
\{1,\ldots,\ell(\sigma_{P})\}$. Since $S(c,|c|,0)=m^{m-2}$ for every cycle $c$ of size $m$, we find
the formula
 $$S((1\ldots n),k,0)=\sum_{(s_{1},\ldots,s_{n})} \frac{n!}{(k+1)!s_{1}!\ldots s_{n}!}  \frac{k!}{1!^{s_{2}}\ldots
(n-1)!^{s_{n}}} 2^{0 s_{2}} 3^{1 s_{3}}\ldots n^{(n-2)s_{n}} ,$$
where the sum is extended to all possible types of partitions of $\{1,\ldots,n\}$ into $n-k$
classes. Let us enumerate the possible types by enumerating the partitions themselves. If we do
this, each type will appear as many times as the number of partitions of this specific type. The
number of partitions of type $(s_{1},\ldots,s_{n})$ is 
$$\frac{n!}{1!^{s_{1}}s_{1}!\ldots n!^{s_{n}}s_{n}!}.$$
Hence, we find the following expression:
$$S((1\ldots n),k,0)= \frac{1}{k+1}\sum_{P} 1^{(1-1)s_{1}} 2^{(2-1)s_{2}}\ldots n^{(n-1)s_{n}},$$
where the sum runs over all partitions of $\{1,\ldots,n\}$ into $n-k$
classes. Now the right hand side has been computed by Kreweras in \cite{Kreweras_CRAS} and it is
equal to
$$S((1\ldots n),k,0)=\frac{1}{k+1}\binom{n-1}{n-k-1}n^k=\binom{n}{k+1}n^{k-1}.$$
This is the expected result. \qed

Let us state separately the following result which has been used in the course of this proof.

\begin{lemma} \label{count paths} Let $P\in NC(n)$ be a partition of type $(s_{1},\ldots,s_{n})$.
There are exactly
$$ \frac{(s_{2}+\ldots+(n-1)s_{n})!}{1!^{s_{2}}\ldots (n-1)!^{s_{n}}} \; 2^{0s_{2}}
\ldots n^{(n-2)s_{n}}$$
increasing paths from $0_{n}$ to $P$.
\end{lemma}

\begin{remark} The explicit expression of the large $N$ limit of the moments of $B_{\frac{t}{N}}$ obtained in this section has been stated by Singer in \cite{Singer} and proved by Biane in \cite{Biane}.

The asymptotic distribution of $B_{\frac{t}{N}}$ as $N$ tends to infinity is the
unique probability measure $\mu_{t}$ on the group $\UC$ of complex numbers of modulus 1 such that,
for all $n\in \N$, $\int_{\UC}\ z^{n} \; \mu(dz)=\int_{\UC} z^{-n} \;
\mu_{t}(dz)=e^{-\frac{nt}{2}}\sum_{k=0}^{n-1} \frac{(-nt)^k}{n k!} \binom{n}{k+1}$. This
expression of the moments of $\mu_{t}$ is not very easy to handle, if only numerically, because it
is an alternated sum of large numbers. Let us point out two analytical ways of studying $\mu_{t}$.

The $S$-transform of $\mu_{t}$ is its moment generating function, defined by $S_{t}(z)=\int_{\UC}
\frac{\xi z}{1-\xi z}\; d\mu_{t}(\xi) = \sum_{n=1}^{\infty} m_{n,t} z^n,$
where $m_{n,t}$ is the $n$-th moment of $\mu_{t}$. Since $|m_{n,t}|\leq 1$ for all $n$ and $t$, the
function $S_{t}$ is holomorphic on the unit disk $\D=\{z: |z|<1\}$. Moreover, since $S_{t}(z)=z+O(z^2)$,
there exist a reciprocal function to $S_{t}$ in a neighbourhood of $0$, which we
denote by $\chi_{t}$. It turns out that $\chi_{t}$ is much simpler than $S_{t}$: 
$\chi_{t}(z)=\frac{z}{z+1}e^{t(z+\frac{1}{2})},$
as one can easily check by using Lagrange's inversion formula. To put it more concisely, the measure $\mu_t$ is fully characterized by the following relation, valid for $z$ in a neighbourhood of $0$:
$$\int_\UC \frac{1}{1-\frac{z}{z+1} e^{tz} e^{\frac{t}{2}}\xi }\; d\mu_t(\xi)=1+z.$$
By studying $\chi_{t}$, Biane
proved in \cite{BianeSB} the following facts. For each $t>0$, the measure $\mu_{t}$ has a
continuous density with
respect to the uniform measure on $\UC$. For $t\in (0,4]$, this density is zero exactly on the set 
$$\left\{e^{i\theta} : |\theta -\pi| \leq 2 \arctan
\sqrt{\frac{4-t}{t}}-\frac{1}{2}\sqrt{t(4-t)}\right\}.$$
For $t>4$, the density of $\mu_{t}$ is positive. Finally, for all $t>0$, the density of $\mu_{t}$
at $e^{i\theta}$ is a real analytic function of $\theta$ on the relative interior of its support. 

Another way of studying $\mu_{t}$ is to observe that $m_{n,t}=\frac{e^{-\frac{nt}{2}}}{2in\pi}
\oint e^{-ntz} \left(1+\frac{1}{z}\right)^n \; dz$, the integral being along any contour of index $1$
with respect to $0$. In \cite{GrossMatytsin}, Gross and Matytsin use this expression and the
saddle point method to exhibit the following phase transition with respect to $t$: if $t\in (0,4)$,
$m_{n,t}$ decays with $n$ like $n^{-\frac{3}{2}}$, whereas for $t\in (4,+\infty)$, it decays
exponentially. Along the same lines, one can check that for $t=4$, $m_{n,4}$ decays like $n^{-\frac{4}{3}}$.
This indicates that the density of $\mu_{4}$ is less regular than the density of $\mu_{t}$ for
$t\in (0,4)$. This behaviour is consistent at a heuristical level with a general result of Biane about additive convolution
with the semi-circle law \cite{BianeConvolution}.
\end{remark}

\subsection{Almost sure convergence} The material gathered so far allows us to prove very easily
the following result.

\begin{proposition} Consider $\sigma \in \S_{n}$. Then, uniformly in $t$ on bounded intervals,
$$\Var\left[p_{\sigma}(B_{\frac{t}{N}})\right]=O(N^{-2}).$$
In particular, on any probability space on which a Brownian motion on $\Un$ is defined for $N$
large enough, the following convergence holds almost surely and in $L^2$ :
$$\lim_{N\to \infty}p_{\sigma}(B_{\frac{t}{N}})=e^{-\frac{nt}{2}}\sum_{k=0}^{n-\ell(\sigma)}
(-1)^k S(\sigma,k,0)\frac{t^k}{k!}.$$
\end{proposition}

\pf Let us denote by $\sigma\times \sigma$ the element of $\S_{2n}$ which sends $i$ on $\sigma(i)$
and $n+i$ on $n+\sigma(i)$ for each $i\in \{1,\ldots,n\}$. With this notation,
$p_{\sigma}^2=p_{\sigma\times \sigma}$, so that
\begin{eqnarray*}
\Var\left[p_{\sigma}(B_{\frac{t}{N}})\right]&=&\E\left[p_{\sigma\times\sigma}(B_{\frac{t}{N}})\right]-
\E\left[p_{\sigma}(B_{\frac{t}{N}})\right]^2\\
&&\hskip -3.5cm =\sum_{k,d=0}^{\infty}
\frac{(-1)^kt^k}{k!N^{2d}}\left[S(\sigma\times\sigma,k,d)-\sum_{\substack{l_{1}+l_{2}=k, 
d_{1}+d_{2}=d}}\frac{k!}{l_{1}!l_{2}!} S(\sigma,l_{1},d_{1})S(\sigma,l_{2},d_{2})\right],
\end{eqnarray*}
where the last expression follows after simplification from Theorem \ref{main expansion}. In the
term corresponding to $d=0$, both $d_{1}$ and $d_{2}$ must be equal to $0$. By the same argument
of support as in the proof of Proposition \ref{fact S}, we find
$$ S(\sigma\times\sigma,k,0)=\sum_{l_{1}+l_{2}=k} \frac{k!}{l_{1}!l_{2}!}
S(\sigma,l_{1},0)S(\sigma,l_{2},0).$$
The result follows immediately.\qed

\section{Asymptotic freeness}

In this section, we prove that independent Brownian motions on $\Un$ converge in distribution, as
$N$ tends to infinity, towards free non-commutative random variables. We do not consider
$*$-freeness, that is, we do not consider products involving $B_{\frac{t}{N}}^{-1}$. In fact,
the asymptotic $*$-freeness of independent Brownian motions follows from a general result of
Voiculescu (see \cite{Voiculescu} and \cite[Lemma 6]{Biane} for details). Here we use Speicher's
characterization of freeness by the vanishing of mixed free cumulants. This
combinatorial characterization is very well suited to the approach we have adopted in this paper.

\subsection{The factorization property}

Let us start by slightly improving Proposition \ref{factorization}, by including the extra
deterministic matrices $M_{1},\ldots,M_{n}$ of Theorem \ref{main expansion}. We consider these
matrices as elements of the non-commutative probability space $(\M_{N}(\C),\tr_{N})$ and speak of
their distribution accordingly.

\begin{proposition}
Let $(M_{1}^{(N)},\ldots,M_{n}^{(N)})_{N\geq 1}$ be a sequence of families of $N\times N$ matrices.
Assume that this sequence converges in distribution. Consider $\sigma \in \S_{n}$. Write $\sigma$
as a product of cycles : $\sigma=(i_{1,1}\ldots i_{1,m_{1}})\ldots (i_{\ell(\sigma),1}\ldots
i_{\ell(\sigma),m_{\ell(\sigma)}})$.
Then
\begin{equation}\label{fact M}
\lim_{N\to \infty}
\E\left[p_{\sigma}(M_{1}^{(N)}B_{\frac{t}{N}},\ldots,M_{n}^{(N)}B_{\frac{t}{N}})\right] =
\prod_{r=1}^{\ell(\sigma)} \lim_{N\to\infty} \E\;\tr_{N}\left(M_{i_{r,1}}^{(N)}B_{\frac{t}{N}}\ldots
M_{i_{r,m_{r}}}^{(N)}B_{\frac{t}{N}}\right).
\end{equation}
More precisely,
\begin{equation}\label{fact M O}
\E\left[p_{\sigma}(M_{1}^{(N)}B_{\frac{t}{N}},\ldots,M_{n}^{(N)}B_{\frac{t}{N}})\right] -
\prod_{r=1}^{\ell(\sigma)} \E\;\tr_{N}\left(M_{i_{r,1}}^{(N)}B_{\frac{t}{N}}\ldots
M_{i_{r,m_{r}}}^{(N)}B_{\frac{t}{N}}\right) = O(N^{-2}),
\end{equation}
uniformly in $t$ on bounded intervals. It is understood that all limits exist.
\end{proposition}

\pf According to Theorem \ref{main expansion}, the left hand side of (\ref{fact M}) is equal to
\begin{equation}\label{calc 1}
e^{-\frac{nt}{2}}\sum_{k=0}^{\infty} \frac{(-t)^{k}}{k!} \sum_{|\sigma '|=|\sigma|-k} \#
\Pi_{k}(\sigma\to\sigma') \lim_{N\to\infty} p_{\sigma'}(M_{1}^{(N)},\ldots,M_{n}^{(N)}).
\end{equation}
This last limit exists for all $\sigma'$ by the assumption that the family
$(M_{1}^{(N)},\ldots,M_{n}^{(N)})$ converges in distribution. All permutations $\sigma'$ which
contribute to the sum satisfy on one hand $|\sigma'|=|\sigma|-k$, hence $|\sigma' \sigma^{-1}|\geq
k$, and $\# \Pi_{k}(\sigma\to\sigma')>0$, hence $|\sigma' \sigma^{-1}|\leq k$. Hence, only
permutations $\sigma'$ such that $\sigma'\preccurlyeq \sigma$ contribute and (\ref{calc 1}) can be
rewritten as
\begin{equation}\label{calc 2}
e^{-\frac{nt}{2}} \sum_{\sigma'\preccurlyeq \sigma} \frac{(-t)^{|\sigma'\sigma^{-1}|}}{|\sigma'
\sigma^{-1}|!}\#\Pi_{|\sigma' \sigma^{-1}|}(\sigma \to \sigma ') \lim_{N\to\infty}
p_{\sigma}(M_{1}^{(N)},\ldots, M_{n}^{(N)}).
\end{equation}
Let $c_{1},\ldots,c_{\ell(\sigma)}$ denote the cycles of $\sigma$. It is not difficult to check
that the interval $[\id,\sigma]$ in the poset $\S_{n}$ is isomorphic to the product of intervals
$\prod_{r=1}^{\ell(\sigma)}[\id, c_{r}]$ by the mapping
$(\alpha_{1},\ldots,\alpha_{\ell(\sigma)})\mapsto \alpha_{1}\ldots \alpha_{\ell(\sigma)}$.
Consider $\sigma'\preccurlyeq \sigma$ and write $\sigma'=\alpha_{1}\ldots \alpha_{\ell(\sigma)}$
accordingly. Then
$|\sigma'\sigma^{-1}|=|\alpha_{1}c_{1}^{-1}|+\ldots+|\alpha_{\ell(\sigma)}c_{\ell(\sigma)}^{-1}|$.
Moreover, by the same argument of shuffling used in the computation of $S((1\ldots n),k,0)$,
$$\#\Pi_{|\sigma'\sigma^{-1}|}(\sigma \to
\sigma')=\frac{|\sigma'\sigma^{-1}|!}{\prod_{r=1}^{\ell(\sigma)}
|\alpha_{r}c_{r}^{-1}|!}\prod_{r=1}^{\ell(\sigma)} \#\Pi_{|\alpha_{r}c_{r}^{-1}|}(c_{r}\to
\alpha_{r}).$$
Hence, (\ref{calc 2}) can be written as
\begin{equation}
\prod_{r=1}^{\ell(\sigma)} e^{-\frac{(|c_{r}|+1)t}{2}} \sum_{\alpha_{r}\preccurlyeq c_{r}}
\frac{(-t)^{|\alpha_{r}c_{r}^{-1}|}}{|\alpha_{r}c_{r}^{-1}|!}
\#\Pi_{|\alpha_{r}c_{r}^{-1}|}(c_{r}\to \alpha_{r}) \lim_{N\to \infty}\tr_{N}(M_{i_{r,1}}^{(N)}\ldots
M_{i_{r,m_{r}}}^{(N)}),
\end{equation}
and this is just the right-hand side of (\ref{fact M}).

The equation (\ref{fact M O}) follows from (\ref{fact M}) just as in Proposition \ref{factorization}.\qed

\subsection{Free cumulants} As a preliminary to the proof of the asymptotic freeness, we compute
the free cumulants of the limiting distribution of $B_{\frac{t}{N}}$. Let $({\mathcal{A}},\varphi)$
be a non-commutative probability
space and $u_{t}$ an element of ${\mathcal{A}}$ such that $B_{\frac{t}{N}}$ converges in
distribution to $u_{t}$ as $N$ tends to infinity. We have spent a substantial part of this paper
proving that the moments of $u_{t}$ are given by
\begin{equation}\label{moments}
\varphi(u_{t}^{n})=e^{-\frac{nt}{2}}\sum_{k=0}^{n-1} \binom{n}{k+1} \frac{(-nt)^{k}}{n k!}.
\end{equation}
Given $a\in {\mathcal A}$ and $\sigma\in \S_{n}$ with cycle lengths $(m_{1},\ldots,m_{\ell(\sigma)})$,
let us use the notation
$\varphi_{\sigma}(a)=\varphi(a^{m_{1}})\ldots \varphi(a^{m_{\ell(\sigma)}})$. The free cumulants
of $u_{t}$ form a family of complex numbers $(k_{\pi}(u_{t}))_{\pi\in \bigcup_{n\geq 1}\S_{n}}$ and
they are characterized by the identity
\begin{equation}\label{mom cum}
\forall n\geq 1 , \forall\sigma \in \S_{n} \; , \;\;
\varphi_{\sigma}(u_{t})=\sum_{\sigma'\preccurlyeq \sigma} k_{\sigma'}(u_{t}).
\end{equation}
Let us use the notation $k_{n}=k_{(1\ldots n)}$. It is an elementary property of the free
cumulants that they are multiplicative, in that $k_{\sigma}=k_{m_{1}}\ldots k_{m_{\ell(\sigma)}}$
when $(m_{1},\ldots,m_{\ell(\sigma)})$ are the cycle lengths of $\sigma$. 

\begin{proposition} The free cumulants of $u_{t}$ are given by
\begin{equation}\label{cumulants}
k_{n}(u_{t})=e^{-\frac{nt}{2}}\frac{(-nt)^{n-1}}{n(n-1)!}.
\end{equation}
More generally, if $\sigma\in \S_{n}$, then
\begin{equation}\label{cumulants pi}
k_{\sigma}(u_{t})=e^{-\frac{nt}{2}}\frac{(-t)^{|\sigma|}}{|\sigma|!}\#\Pi_{|\sigma|}(\id
\to\sigma).
\end{equation} 
\end{proposition}

\pf Let us put $\varphi(u_{t}^{n})$ under the form  of the right hand side of (\ref{mom cum}).
Applying Theorem \ref{main expansion}, using Kreweras
complementation and using Lemma \ref{count paths}, we find
\begin{eqnarray*}
\varphi(u_{t}^{n})&=&e^{-\frac{nt}{2}}\sum_{\sigma \preccurlyeq (1\ldots n)} \frac{(-t)^{|\sigma
(1\ldots n)^{-1}|}}{|\sigma (1\ldots n)^{-1}|!}\#\Pi_{|\sigma (1\ldots n)^{-1}|}((1\ldots n)\to
\sigma) \\
&=& e^{-\frac{nt}{2}}\sum_{\sigma \preccurlyeq (1\ldots n)}
\frac{(-t)^{|\sigma|}}{|\sigma|!}\#\Pi_{|\sigma|}(\id \to \sigma) \\
&=&\sum_{\sigma \preccurlyeq (1\ldots n)} \frac{(e^{-\frac{1t}{2}}(-t)^{0}1^{-1})^{s_{1}}(e^{-\frac{2t}{2}}(-t)^{1}2^{0})^
{s_{2}}\ldots (e^{-\frac{nt}{2}}(-t)^{n-1}n^{n-2})^{s_{n}}}{0!^{s_{1}}1!^{s_{2}}\ldots
(n-1)!^{s_{n}}},
\end{eqnarray*}
where $s_{1},\ldots,s_{n}$ are respectively the number of fixed points of $\sigma$, and the
numbers of transpositions, 3-cycles, $\ldots$, $n$-cycles in the decomposition of $\sigma$.
By comparing this expression with (\ref{mom cum}), we find the desired expression for the
cumulants of $u_{t}$. \qed

Let us recall briefly Speicher's characterization of freeness by the vanishing of mixed free
cumulants \cite{Speicher}. Let $a_{1},\ldots,a_{n}$ be non-commutative random variables on a space
$(\mathcal{A},\varphi)$, where $\varphi$ is a tracial state.
For $\sigma\in \S_{n}$, the number $m_{\sigma}(a_{1},\ldots,a_{n})$ is defined by
$$m_{\sigma}(a_{1},\ldots,a_{n})=\prod_{\substack{c \mbox{ \scriptsize  cycle of } \sigma\\ c=(i_{1}\ldots
i_{r})}} \varphi(a_{i_{1}}\ldots a_{i_{r}}).$$
It is well defined thanks to the fact that $\varphi$ is tracial. The numbers
$m_{\sigma}(a_{1},\ldots,a_{n})$ are the mixed moments of $a_{1},\ldots,a_{n}$. The relation
$$m_{\sigma}(a_{1},\ldots,a_{n})=\sum_{\sigma'\preccurlyeq
\sigma}k_{\sigma'}(a_{1},\ldots,a_{n})$$
characterizes the family of numbers $k_{\sigma}(a_{1},\ldots,a_{n})$, the mixed free cumulants of
$a_{1},\ldots,a_{n}$. 

Speicher's characterization of freeness is the following. Let $(a_{k})_{k\geq 1}$ be a family of
elements of ${\mathcal A}$. Then this family is free if and only if, for all $n\geq 2$ and all
$i_{1},\ldots,i_{n}\geq 1$ such that $i_{r}\neq i_{s}$ for some $r,s\in \{1,\ldots,n\}$,
$k_{(1\ldots n)}(a_{i_{1}},\ldots,a_{i_{n}})=0$.

\begin{theorem} Let $(B^{(N,k)})_{N,k\geq 1}$ be a family of Brownian motions, such that, for all
$N,k\geq 1$, $B^{(N,k)}$ is a Brownian motion on $\Un$ and, for all $N\geq 1$, the Brownian
motions $(B^{(N,k)})_{k\geq 1}$ are independent. Let $(t_{k})_{k\geq 1}$ be a sequence of
non-negative real numbers.

Then, as $N$ tends to infinity, the family of non-commutative random variables
$(B^{(N,k)}_{\frac{t_{k}}{N}})_{k\geq 1}$ converges in distribution towards a free family
$(b^{(k)})_{k\geq 1}$ of non-commutative random variables such that, for all $k\geq 1$, $b^{(k)}$
has the distribution of $u_{t_{k}}$.
\end{theorem}

\pf We prove the result for finite families $(B^{(N,k)})_{N\geq 1, k\in
\{1,\ldots,K\}}$ for some finite $K$, by induction on $K$. The case $K=1$ is settled by our
computation of the asymptotic distribution of $B^{(N)}_{\frac{t}{N}}$, that is, Propositions \ref{asympt
dist} and \ref{calc S}.

Let $K\geq 2$ be an integer and let us assume that the property is proved for $K-1$ independent
Brownian motions. We need to prove that the mixed free cumulants of
$B^{(N,1)}_{\frac{t_{1}}{N}},\ldots,B^{(N,K)}_{\frac{t_{K}}{N}}$ tend to zero as $N$ tends to
infinity. We regard $B^{(N,k)}_{\frac{t_{k}}{N}}$ as elements of the non-commutative probability 
space $(L^{\infty}(\Omega,\P)\otimes\M_{N}(\C),\E \otimes \tr_{N})$ and we use the notation
$m_{\sigma}$ and $k_{\sigma}$ accordingly. In particular, with our previous notation, $m_{\sigma}=\E
\; p_{\sigma}$.

What we need to prove is that, for all $n\geq 2$, all $\sigma \in \S_{n}$, all $i_{1},\ldots, i_{n}\in
\{1,\ldots,K\}$ not all equal, 
$$ \lim_{N\to\infty} m_{\sigma}(B^{(N,i_{1})}_{\frac{t_{i_{1}}}{N}},\ldots,B^{(N,i_{n})}_{\frac{t_{i_{n}}}{N}}) =
\sum_{\substack{\sigma'\preccurlyeq \sigma \\ \forall r\in\{1,\ldots,n\},\; i_{\sigma'(r)}=i_{r}}}
\lim_{N\to\infty}
k_{\sigma'}(B^{(N,i_{1})}_{\frac{t_{i_{1}}}{N}},\ldots,B^{(N,i_{n})}_{\frac{t_{i_{n}}}{N}}).$$
By the factorization property (\ref{fact M}) and Fubini's theorem, the left hand side is
multiplicative with respect to the cycle decomposition of $\sigma$. The right hand side is also
clearly multiplicative, hence, it suffices to consider the case where $\sigma=(1\ldots n)$. In
this case, we are looking at the expected trace of a product $B^{(N,i_{1})}_{\frac{t_{i_{1}}}{N}}\ldots
B^{(N,i_{n})}_{\frac{t_{i_{n}}}{N}}$ where at least two factors are distinct. Of course, the case
where one of the $K$ possible factors does not appear is treated by the induction hypothesis. Let
us assume that the $K$ factors appear, in particular $B^{(N,1)}$. Up to cyclic permutation, which
does not affect its trace, the product above can be put under the form $W_{1}^{(N)}B^{(N,1)}_{\frac{t_{1}}{N}}\ldots
W_{r}^{(N)}B^{(N,1)}_{\frac{t_{1}}{N}}$, for some $r\geq 1$ and some products
$W_{1}^{(N)},\ldots,W_{r}^{(N)}$ of factors among $B^{(N,2)},\ldots,B^{(N,K)}$. Our previous
results show that
\begin{eqnarray*}
\lim_{N\to\infty} \E \;\tr_{N}\left(B^{(N,i_{1})}_{\frac{t_{i_{1}}}{N}}\ldots
B^{(N,i_{n})}_{\frac{t_{i_{n}}}{N}}\right)&=& e^{-\frac{nt}{2}}\sum_{\sigma\preccurlyeq (1\ldots r)}
\frac{(-t)^{|\sigma (1\ldots r)^{-1}|}}{|\sigma (1\ldots r)^{-1}|!}\\
&&\hskip -0cm  \#\Pi_{|\sigma (1\ldots n)^{-1}|}((1\ldots r)\to \sigma) \lim_{N\to
\infty} m_{\sigma}(W_{1}^{(N)},\ldots,W_{r}^{(N)})\\
&&\hskip -2cm = e^{-\frac{nt}{2}}\sum_{\sigma\preccurlyeq (1\ldots r)}
\frac{(-t)^{|\sigma |}}{|\sigma |!} \#\Pi_{|\sigma|}(\id \to \sigma)
\lim_{N\to\infty} m_{K(\sigma)}(W_{1}^{(N)},\ldots,W_{r}^{(N)})\\
&&\hskip -2cm = \sum_{\sigma\preccurlyeq (1\ldots r)}
k_{\sigma}(u_{t_{1}})\lim_{N\to\infty}m_{K(\sigma)} (W_{1}^{(N)},\ldots,W_{r}^{(N)}),
\end{eqnarray*}
where we have changed $\sigma$ in $K(\sigma)=\sigma^{-1}(1\ldots r)$ between the first and the
second line. 

The term $\lim_{N\to\infty}m_{K(\sigma)} (W_{1}^{(N)},\ldots,W_{r}^{(N)})$ is equal to a
sum of limits of free cumulants of the factors appearing in $W_{1}^{(N)},\ldots,W_{r}^{(N)}$ in
this order. By induction, only pure free cumulants appear, those which do not involve more than
one Brownian motion in each cycle of the permutation. Moreover, by the combinatorial description 
of $K(\sigma)$ as a non-crossing partition, only such cumulants appear that remain
non-crossing when they are merged with $\sigma$. In symbols,
$$\lim_{N\to\infty} \E \;\tr_{N}\left(B^{(N,i_{1})}_{\frac{t_{i_{1}}}{N}}\ldots
B^{(N,i_{n})}_{\frac{t_{i_{n}}}{N}}\right)=\sum_{\substack{\sigma \preccurlyeq (1\ldots n)\\ \forall
k\in \{1,\ldots,n\},\; i_{k}=i_{\sigma(k)}}} \lim_{N\to\infty}
k_{\sigma}(B^{(N,i_{1})}_{\frac{t_{i_{1}}}{N}},\ldots,
B^{(N,i_{n})}_{\frac{t_{i_{n}}}{N}}).$$
This is exactly what we expected. \qed

\section{Large N Yang-Mills theory on a disk and branching covers}

In this section, we explain how our main expansion relates the Brownian motion on the unitary group
to a natural model of random branching covers on this disk. In doing this, we give a rigorous
proof of results which are stated in \cite{GrossTaylor2}.

Let $D$ be the closed disk of radius $1$ centred at the origin $O$ of $\R^2$. Let $n\geq 1$ be an
integer. Let $\lambda$ be a partition of $n$. We define the set $\Ram_{n,\lambda}$ as the set of
isomorphism classes of ramified coverings $\pi:R \lra D$ which satisfy the following properties.\\
1. $R$ is a ramified covering of degree $n$.\\
2. For each ramification point $x\in D$ of $R$ which is not the origin $O$, $R$ has a generic
ramification type at $x$, in that $\#\pi^{-1}(x)=n-1$.\\
3. The monodromy of $R$ along the boundary of $D$ belongs to the conjugacy class of $\S_{n}$
corresponding to $\lambda$.\\
An element of $\Ram_{n,\lambda}$ is allowed to be ramified over $O$, with any kind of
ramification. The set of its ramification points distinct from $O$ is called its locus of generic
ramification. It is contained in the interior of $D-\{0\}$, which we denote by $D^*$.

Let $X$ be a finite subset of $D^*$. We define $\Ram_{n,\lambda,X}$ as the subset of
$\Ram_{n,\lambda}$ formed by the coverings whose locus of generic ramification is $X$.

The set $\Ram_{n,\lambda,X}$ is in natural one-to-one correspondence with a set of equivalence
classes of paths in the Cayley graph of $\S_{n}$ as follows. Assume
that $X=\{x_{1},\ldots,x_{k}\}$. Choose a point $b$ on the boundary of $D$. By the interior of a
simple closed continuous curve based at $b$, we mean the bounded connected component of the
complement of its range. Let $C,C_{1},\ldots,C_{k}$ be simple closed curves in $D$ based at $b$
with pairwise disjoint interiors such that the interior of $C$ contains $O$ and, for each $i\in\{1,\ldots,k\}$, the
interior of $C_{i}$ contains $x_{i}$. We assume that this is done in such a way that the curve
$CC_{k}\ldots C_{1}$ is homotopic to the boundary of $D$ in $D-(X\cup\{O\})$. Then the monodromies
$\sigma_{O},\tau_{1},\ldots,\tau_{k}$ of an element $R\in \Ram_{n,\lambda,X}$ along the curves $C,C_{1},\ldots,C_{k}$
are defined in $\S_{n}$ up to simultaneous conjugation and their orbit characterizes $R$. The
assumptions made on $R$ imply that $\tau_{1},\ldots,\tau_{k}$ are transpositions and
$\sigma_{O}\tau_{1}\ldots \tau_{k}$ belongs to $\lambda$.

Let $\PP_{n,\lambda,k}$ be the set of paths in the Cayley graph of $\S_{n}$ which start at an
element of the conjugacy class determined by $\lambda$. The symmetric group acts on
$\PP_{n,\lambda,k}$ by conjugation. The mapping from $\Ram_{n,\lambda,X}$ which associates to $R$
the orbit of the path $(\sigma_{O}\tau_{k}\ldots \tau_{1},\sigma_{O}\tau_{k}\ldots \tau_{2}, 
\ldots,\sigma_{O}\tau_{k},\sigma_{O})$ is a bijection. Moreover, the cardinal of the stabilizer of
this orbit is equal to the cardinal of the automorphism group $\Aut(R)$ of $R$. Hence, the image on
$\Ram_{n,\lambda,X}$ of the counting measure on $\PP_{n,\lambda,k}$ by the mapping
 $\PP_{n,\lambda,k}\lra \PP_{n,\lambda,k}/\S_{n}\simeq \Ram_{n,\lambda,X}$ is the measure
$$\rho_{n,\lambda,X}=\sum_{R\in \Ram_{n,\lambda}(X)} \frac{n!}{\#\Aut(R)} \; \delta_{R}.$$
This measure is finite and satisfies $\rho_{n,\lambda,X}(1)=\binom{n}{2}^k$.

There is a natural topology on $\Ram_{n,\lambda}$, which is generated by the sets
$$\mathcal{O}(R,U)=\{R'\in \Ram_{n,\lambda} : R_{|D-U}\simeq R'_{|D-U}\},$$
where $R$ spans $\Ram_{n,\lambda}$ and $U$ the set of neighbourhoods of the locus of generic ramification
of $R$. This is a fairly coarse topology: for example, the cardinal of the locus of generic ramification
is not continuous, but only lower semi-continuous in this topology. However, the ramification
type at $O$ is continuous. On the set ${\mathcal X}$ of
finite subsets of $D^*$, we put the topology which makes the bijection ${\mathcal X}\simeq \bigsqcup_{n\geq 0}
((D^*)^n-\Delta_{n})/\S_{n}$ a homeomorphism, where $\Delta_{n}$ is the subset of $(D^*)^n$ where at
least two components coincide. These topologies do not make the ramification locus a continuous
function of the ramified covering. Nevertheless, let ${\mathcal M}(\Ram_{n,\lambda})$ denote the
space of finite Borel measures on $\Ram_{n,\lambda}$ endowed with the topology of weak convergence.

\begin{lemma} The mapping ${\mathcal X}\lra {\mathcal M}(\Ram_{n,\lambda})$ which sends $X$ to
$\rho_{n,\lambda,X}$ is continuous.
\end{lemma}

\pf By definition of the topology on ${\mathcal X}$, it suffices to prove that the mapping is
continuous on $(D^*)^k-\Delta_{k}$ for all $k\geq 0$. Consider $k\geq 0$, $X=\{x_{1},\ldots,x_{k}\}$
and a bounded continuous function $f:\Ram_{n,\lambda,X}\lra \R$. Choose $\epsilon>0$. 

Since $\Ram_{n,\lambda,X}$ is finite, the continuity of $f$ implies the existence of $r>0$ such
that the balls $B(x_{i},r)$ are contained in $D^*$, pairwise disjoint for $i\in\{1,\ldots,k\}$
and the neighbourhood $U=B(x_{1},r)\times  \ldots \times B(x_{k},r)$ of $X$ in $(D^*)^k-\Delta_{k}$
satisfies 
$$\forall R \in \Ram_{n,\lambda,X}, \forall R' \in {\mathcal O}(R,U) ,\;  |f(R')-f(R)| <
\epsilon \binom{n}{2}^{-k}.$$
Let $X'=\{x'_{1},\ldots,x'_{k}\}$ be an element of $U$. Let $\phi$ be a homeomorphism of $D$ such
that $\phi_{|D-U}=\id_{D-U}$ and $\phi(x_{i})=x'_{i}$ for all $i\in\{1,\ldots,k\}$. For each
ramified covering $\pi:R\lra D$ belonging to $\Ram_{n,\lambda,X}$, the covering $\Phi(R)=(\phi\circ
\pi : R\lra D)$ belongs to $\Ram_{n,\lambda,X'}$. Replacing $\phi$ by its inverse in the
definition of $\Phi:\Ram_{n,\lambda,X}\lra \Ram_{n,\lambda,X'}$ yields the inverse mapping, hence $\Phi$
is a bijection. Moreover, the conjugation by $\phi$ determines an isomorphism
between $\Aut(R)$ and $\Aut(\Phi(R))$. Finally, $R$ and $\Phi(R)$ are isomorphic outside $U$.
Altogether,
$$\left| \rho_{n,\lambda,X'}(f)-\rho_{n,\lambda,X}(f)\right| \leq
\sum_{R\in\Ram_{n,\lambda,X}} \frac{n!}{\#\Aut(R)}\; \left|
f(\Phi(R))-f(R)\right|<\epsilon.$$
Since $k$, $X$, $f$ and $\epsilon$ were arbitrary, the result follows. \qed

Let $t\geq 0$ be a real number. Let $\Xi_{t}$ be the distribution of a Poisson point process on $D$
of intensity $\frac{t}{\pi}$ times the Lebesgue measure on $D$. Under $\Xi_{t}$, a random subset of
$D$ is contained in $D^*$ with probability $1$ and the average number of points of such a random
set is $t$. Thinking of $\Xi_{t}$ as a Borel probability measure
on ${\mathcal X}$, we define a measure on $\Ram_{n,\lambda}$ by setting
$$\rho_{n,\lambda}^{t}=\int_{{\mathcal X}} \rho_{n,\lambda,X} \; \Xi_{t}(dX).$$
The measure $\rho_{n,\lambda}^{t}$ is finite and satisfies $\rho_{n,\lambda}^{t}(1)=
e^{t\binom{n}{2}-t}$.
We define a probability measure $\mu_{n,\lambda}^{t}$ on $\Ram_{n,\lambda}$ by normalizing
$\rho_{n,\lambda}^{t}$.

Let us define two functions on $\Ram_{n,\lambda}$. Firstly, given $R\in \Ram_{n,\lambda}$, let us
define $k(R)$ as the number of ramification points of $R$ distinct from $O$. We have observed that
this is a lower semi-continuous, hence measurable function of $R$. Secondly, let $\chi(R)$ be the
Euler characteristic of $R$. 

\begin{lemma} \label{chi} Let $X$ be a subset of cardinal $k$ of $D^*$. Let $R$ be an element
of $\Ram_{n,\lambda,X}$ and $\gamma=(\sigma_{0},\ldots,\sigma_{k})$ a representative of the
associated orbit of $\PP_{n,\lambda,k}$. Then 
$$\chi(R)=\ell(\sigma_{k})-k=\ell(\lambda)-2d(\gamma).$$
In particular, $\chi:\Ram_{n,\lambda}\lra \Z$ is upper semi-continuous and measurable.
\end{lemma} 

\pf The first equality follows from the Riemann-Hurwitz formula, the second from the definition of
the defect of a path. The last assertion follows from the lower
semi-continuity of $k$ and the fact that $\ell(\sigma_{k})$, which
depends only on the ramification type at $O$, is a continuous function of $R$.\qed

The main result is the following.

\begin{theorem} Let $N,n\geq 1$ be two integers. Let $(B_{t})_{t\geq 0}$ be the Brownian motion on
$U(N)$ defined in Theorem \ref{main expansion}. Let $\lambda$ be a partition of $n$ and $\sigma$
an element of $\S_{n}$ which belongs to the conjugacy class determined by $\lambda$.  Let $t\geq
0$ be a real number. Then
$$e^{nt-\frac{n^2t}{2}}\E\left[p_{\sigma}^{st}(B_{\frac{t}{N}})\right]=
\int_{\Ram_{n,\lambda}} (-1)^{k(R)} N^{\chi(R)}\; \mu_{n,\lambda}^{t}(dR).$$
\end{theorem}

We could have avoided the unpleasant exponential factor in the statement of this theorem if we had
considered the signed measure $\tilde \rho_{n,\lambda}^{t}=\int_{{\mathcal X}} (-1)^{\# X} \rho_{n,\lambda,X} \;
\Xi_{t}(dX)$ instead of $\rho_{n,\lambda}$. \\

\pf Let $X$ be a finite subset of $D^*$. The set $\Ram_{n,\lambda,X}$ is in
bijection with the set of orbits of $\Pi_{k}(\sigma)$ under the action of $\S_{n}$ by conjugation.
Let $R$ be an element of $\Ram_{n,\lambda,X}$ and $\gamma=(\sigma,\sigma_{1},\ldots,\sigma_{k})$ a
representative of the corresponding orbit. By Lemma \ref{chi},
\begin{eqnarray*}
\int_{\Ram_{n,\lambda}} (-1)^{k(R)} N^{\chi(R)} \; \rho_{n,\lambda,X}(dR) &=& \sum_{\gamma\in
\Pi_{k}(\sigma)} (-1)^{k(X)} N^{\ell(\lambda)-2d(\gamma)}\\
&=& N^{\ell(\lambda)} \sum_{d\geq 0} \frac{(-1)^{k(X)}}{N^{2d}} S(\sigma,k(X),d).
\end{eqnarray*}
Integrating with respect to $\Xi_{t}(dX)$, we find
$$\int_{\Ram_{n,\lambda}} (-1)^{k(R)} N^{\chi(R)} \; \rho^{t}_{n,\lambda}(dR) =
e^{-t}N^{\ell(\sigma)} \sum_{k,d\geq 0} \frac{(-t)^{k}}{k!N^{2d}} S(\sigma,k,d).$$
By Theorem \ref{main expansion}, the right-hand side of this equality is equal to
$e^{-t+\frac{nt}{2}}\E\left[p_{\sigma}^{st}(B_{\frac{t}{N}})\right]$. The result follows after
normalizing $\rho_{n,\lambda}^{t}$. \qed

\bibliographystyle{elsart-num}
\bibliography{ymsw}

\begin{thebibliography}{10}
\expandafter\ifx\csname url\endcsname\relax
  \def\url#1{\texttt{#1}}\fi
\expandafter\ifx\csname urlprefix\endcsname\relax\def\urlprefix{URL }\fi

\bibitem{Biane}
P.~Biane, Free {B}rownian motion, free stochastic calculus and random matrices,
  in: Free probability theory (Waterloo, ON, 1995), Vol.~12 of Fields Inst.
  Commun., Amer. Math. Soc., Providence, RI, 1997, pp. 1--19.

\bibitem{Xu}
F.~Xu, A random matrix model from two-dimensional {Y}ang-{M}ills theory, Comm.
  Math. Phys. 190~(2) (1997) 287--307.

\bibitem{Collins}
B.~Collins, Moments and cumulants of polynomial random variables on unitary
  groups, the {I}tzykson-{Z}uber integral, and free probability, Int. Math.
  Res. Not. ~(17) (2003) 953--982.

\bibitem{CollinsSniady}
B.~Collins, P.~{\'S}niady, Integration with respect to the {H}aar measure on
  unitary, orthogonal and symplectic group, Comm. Math. Phys. 264~(3) (2006)
  773--795.

\bibitem{KK}
V.~A. Kazakov, I.~K. Kostov, Nonlinear strings in two-dimensional {${\rm
  U}(\infty )$} gauge theory, Nuclear Phys. B 176~(1) (1980) 199--215.

\bibitem{GrossTaylor2}
D.~J. Gross, W.~Taylor, Twists and {W}ilson loops in the string theory of two
  dimensional {QCD}, Nuclear Physics B 403 (1993) 395.

\bibitem{GrossMatytsin}
D.~J. Gross, A.~Matytsin, Some properties of large-{$N$} two-dimensional
  {Y}ang-{M}ills theory, Nuclear Phys. B 437~(3) (1995) 541--584.

\bibitem{GopakumarGross}
R.~Gopakumar, D.~J. Gross, Mastering the master field, Nuclear Phys. B
  451~(1-2) (1995) 379--415.

\bibitem{Singer}
I.~M. Singer, On the master field in two dimensions, in: Functional analysis on
  the eve of the 21st century, Vol.\ 1 (New Brunswick, NJ, 1993), Vol. 131 of
  Progr. Math., Birkh\"auser Boston, Boston, MA, 1995, pp. 263--281.

\bibitem{Sengupta}
A.~N. Sengupta, Traces in two-dimensional {QCD}: the large-{$N$} limit, in:
  Traces in geometry, number theory and quantum fields (edited by Sergio
  Albeverio, Matilde Marcolli, Sylvie Paycha, and Jorge Plazas), Vieweg (to
  appear).

\bibitem{GoodmanWallach}
R.~Goodman, N.~R. Wallach, Representations and invariants of the classical
  groups, Vol.~68 of Encyclopedia of Mathematics and its Applications,
  Cambridge University Press, Cambridge, 1998.

\bibitem{Macdonald}
I.~G. Macdonald, Symmetric functions and {H}all polynomials, 2nd Edition,
  Oxford Mathematical Monographs, The Clarendon Press Oxford University Press,
  New York, 1995, with contributions by A. Zelevinsky, Oxford Science
  Publications.

\bibitem{Liao}
M.~Liao, L\'evy processes in {L}ie groups, Vol. 162 of Cambridge Tracts in
  Mathematics, Cambridge University Press, Cambridge, 2004.

\bibitem{VershikOkounkov}
A.~Okounkov, A.~Vershik, A new approach to representation theory of symmetric
  groups, Selecta Math. (N.S.) 2~(4) (1996) 581--605.

\bibitem{LandoZvonkin}
S.~K. Lando, A.~K. Zvonkin, Graphs on surfaces and their applications, Vol. 141
  of Encyclopaedia of Mathematical Sciences, Springer-Verlag, Berlin, 2004,
  with an appendix by Don B. Zagier, Low-Dimensional Topology, II.

\bibitem{OLEIS}
N.~Sloane, The On-Line Encyclopedia of Integer Sequences,
  www.research.att.com/~njas/sequences/, 2007.

\bibitem{Kreweras}
G.~Kreweras, Sur les partitions non crois\'ees d'un cycle, Discrete Math. 1~(4)
  (1972) 333--350.

\bibitem{Kreweras_CRAS}
G.~Kreweras, Une famille d'identit\'es mettant en jeu toutes les partitions
  d'un ensemble fini de variables en un nombre donn\'e de classes, C. R. Acad.
  Sci. Paris S\'er. A-B 270 (1970) A1140--A1143.

\bibitem{BianeSB}
P.~Biane, Segal-{B}argmann transform, functional calculus on matrix spaces and
  the theory of semi-circular and circular systems, J. Funct. Anal. 144~(1)
  (1997) 232--286.

\bibitem{BianeConvolution}
P.~Biane, On the free convolution with a semi-circular distribution, Indiana
  Univ. Math. J. 46~(3) (1997) 705--718.

\bibitem{Voiculescu}
D.~V. Voiculescu, K.~J. Dykema, A.~Nica, Free random variables, Vol.~1 of CRM
  Monograph Series, American Mathematical Society, Providence, RI, 1992, a
  noncommutative probability approach to free products with applications to
  random matrices, operator algebras and harmonic analysis on free groups.

\bibitem{Speicher}
R.~Speicher, Free probability theory and non-crossing partitions, S\'em.
  Lothar. Combin. 39 (1997) Art.\ B39c, 38 pp.\ (electronic).

\end{thebibliography}






\end{document}